\documentclass[aap]{imsart}
\RequirePackage[OT1]{fontenc}
\RequirePackage{amsthm,amsmath}
\usepackage{amsfonts}
\usepackage{amssymb}
\usepackage{color}
\usepackage{epsfig}
\def\inter{{\mathop{\hbox{\rm int}}}}
\def\Dom{{{\mathop{\hbox{\rm Dom}}}}}
\def\argmin{\mathop{\hbox{\rm argmin\,}}}

\def\bR{{\mathbf{R}}}
\def\bS{{\mathbf{S}}}
\def\Tr{\mathop{\hbox{\rm Tr}}}
\def\Diag{\mathop{\hbox{\rm Diag}}}

\def\cF{{\cal F}}

\def\cC{{\cal C}}

\def\cF{{\cal F}}
\def\bS{{\mathbf{S}}}

\def\bE{{\mathbf{E}}}
\def\Prob{{\hbox{\rm Prob}}}

\def\argmin{\mathop{\rm argmin}}

\newcommand{\g} {\gamma}

\newcommand{\te}{\theta}

\newcommand{\sign}{\mbox{sign}}

\newcommand{\be}{\begin{eqnarray}}
\newcommand{\ee}[1]{\label{eq:#1}\end{eqnarray}}
\newcommand{\nn}{\nonumber \\}
\newcommand{\ese}{\end{eqnarray*}}
\newcommand{\bse}{\begin{eqnarray*}}
\newcommand{\rf}[1]{~(\ref{eq:#1})}

\newtheorem{example}{Example}[section]
\newtheorem{theorem}{Theorem}[section]
\newtheorem{lemma}{Lemma}

\newtheorem{proposition}{Proposition}[section]

\newtheorem{definition}{Definition}[section]

\def\blacksquare{\hbox{\vrule width 4pt height 4pt depth 0pt}}

\newtheorem{theo}{Theorem}[section]
\newcommand{\bthm}{\begin{theo}}
\newcommand{\ethm}[1]{\label{the:#1}\end{theo}\par}

\newtheorem{col}{Corollary}[section]
\newcommand{\bcol}{\begin{col}}
\newcommand{\ecol}[1]{\label{the:#1}\end{col}\par}

\newtheorem{defi}{Definition}[section]
\newcommand{\bdf}{\begin{defi}}
\newcommand{\edf}[1]{\label{df:#1}\end{defi}\par}

\newtheorem{lem}{Lemma}[section]
\newcommand{\blem}{\begin{lem}}
\newcommand{\elem}[1]{\label{le:#1}\end{lem}\par}

\newtheorem{pro}{Proposition}[section]
\newcommand{\bpro}{\begin{pro}}
\newcommand{\epro}[1]{\label{pr:#1}\end{pro}\par}

\newtheorem{statement}{Problem}[section]
\newcommand{\bstatement}{\begin{statement}}
\newcommand{\estatement}[1]{\label{stat:#1}\end{statement}\par}

\newcounter{assc}[section]
\newcommand{\bass}[1]{\refstepcounter{assc}\label{ass:#1} \begin{maliste}{\bf \arabic{assc}.}}
\newcommand{\eass}{\end{maliste}}

\newcounter{algoc}
\newcommand{\balgo}[1]{\refstepcounter{algoc}\label{algo:#1} \begin{malistea}{\bf \arabic{algoc}.}}
\newcommand{\ealgo}{\end{malistea}}

\newcommand{\myqed}{\hfill\hbox{\hskip 4pt
                \vrule width 5pt height 4pt depth 1.5pt}\vspace{0.5cm}\par}
\def    \qed    {\hfill\hbox{\hskip 4pt
                \vrule width 3pt height 1pt depth 1pt
                \hbox{\vrule width 2pt height 3.5pt depth 1pt}}}

\newcommand{\aic}[2]{{\color{blue}~#2}}
\endlocaldefs
\begin{document}

\begin{frontmatter}
\title{Large Deviations of Vector-valued Martingales in 2-Smooth Normed Spaces}
\runtitle{Large Deviations of Vector-valued Martingales}
\begin{aug}
\author{\fnms{Anatoli B.} \snm{Juditsky}\ead[label=e1]{juditsky@imag.fr}}
\and
\author{\fnms{Arkadi S.}\snm{Nemirovski}\thanksref{t1}\ead[label=e2]{nemirovs@isye.gatech.edu}}
\thankstext{t1}{
Research of
this author was partly supported by the NSF award DMI-0619977.}
\runauthor{A. Juditsky and A. Nemirovski}
\affiliation{Universit\'e  Grenoble I and Georgia Institute of Technology}
\address{Anatoli Juditsky\\
Laboratoire Jean Kuntzmann \\
51 rue des Math\'ematiques\\
BP 53\\
38041 Grenoble Cedex 9
France \\
E-mail: \printead{e1}}
\address{Arkadi Nemirovski\\
School of Industrial and Systems Engineering\\
Georgia Institute of Technology\\
765 Ferst Drive\\
 Atlanta GA 30332-0205 USA\\
\ E-mail: \printead{e2}}
\end{aug}

\begin{abstract} In this paper, we derive exponential
bounds on probabilities of large deviations for ``light tail''
martingales taking values in finite-dimensional normed
spaces. Our primary emphasis is on the case where the bounds are
dimension-independent or nearly so. We demonstrate that this is
the case when the norm on the space can be approximated, within an
absolute constant factor, by a norm which is differentiable on the
unit sphere with a Lipschitz continuous 
gradient. We also present various examples of spaces possessing the latter
property.
\end{abstract}
\medskip



\end{frontmatter}
\section{Introduction}
It is well-known that for a sequence of independent zero mean
random reals $\{\xi_i\}_{i=1}^\infty$ with light tail
distributions (e.g., such that
$\bE\left\{\exp\{|\xi_i|^\alpha\sigma_i^{-\alpha}\}\right\}\leq
\exp\{1\}$ for certain $\alpha\in[1,2]$ and deterministic
$\sigma_t>0$), a ``typical magnitude'' of the sum
$S_t=\sum_{i=1}^t\xi_i$ is ``at most of order of
$\sqrt{\sum_{i=1}^t\sigma_i^2}\;$'', meaning that
$$\Prob\left\{|S_t|>[1+\gamma]\sqrt{\sum_{i=1}^t\sigma_i^2}\right\}\leq
O(1)\exp\{-O(1)\gamma^\alpha\}$$
for all $\gamma\geq0$; here in
what follows, all $O(1)$ are positive absolute constants. The
question we focus on in this paper is to which extent the above
large deviation bound is preserved when passing from scalar random
variables to independent zero mean random variables taking values
in a normed space $(E,\|\cdot\|)$ of (possibly, large) dimension
$n<\infty$. Now our ``light tail'' condition reads
\begin{equation}\label{lighttail}
\bE\left\{\exp\{\|\xi_i\|^\alpha\sigma_i^{-\alpha}\}\right\}\leq\exp\{1\}
\end{equation}
for some $\alpha\in [1,2]$,
and what we want to get is a bound of the form
$$
\forall
\gamma\geq0:\Prob\left\{\|\sum_{i=1}^t\xi_i\|>[\theta+\gamma]\sqrt{\sum_{i=1}^t\sigma_i^2}\right\}\leq
O(1)\exp\{-O(1)\gamma^\alpha\}\eqno{(*)}
$$
with a ``moderate'' value of the constant $\theta$. It is
immediately seen that our goal is not always attainable. For
instance, let $(E,\|\cdot\|)$ be $\ell_1^n$ (i.e., $\bR^n$ equipped
with the norm $\|x\|_1=\sum_{i=1}^n|x_i|$), and let $\xi_i$ take
values $\pm e_i$ with probability $1/2$, $1\leq i\leq n$, where
$e_i$ are the standard basic orths in $\bR^n$. Then
(\ref{lighttail}) holds true with $\sigma_i=1$, while
$\|S_k\|_1\equiv k$ whenever $k\leq n$. We see that in order for
$(*)$ to be true, $\theta$ should be as large as $O(1)\sqrt{n}$.
On the other hand, with $\theta=O(1)\sqrt{\dim E}$, $(*)$ indeed is
true independently of the norm $\|\cdot\|$ in question (see
Example \ref{wwex:R_n} in Section \ref{sect1}). Our major goal in this paper is to
show that a sufficient condition for $(*)$ to be valid with
certain $\theta$ is {\sl $\theta^2$-regularity} of the space
$(E,\|\cdot\|)$. The latter means, essentially, that $\|\cdot\|$ can be
approximated within an absolute constant factor by a norm
$p(\cdot)$ which is continuously differentiable outside of the
origin and possesses Lipschitz continuous, with the Lipschitz constant
$\theta^2$, derivative on its unit sphere:
\begin{equation}\label{eq555}
p(x)=p(y)=1\Rightarrow p_*(p'(x)-p'(y))\leq\theta^2p(x-y)
\end{equation}
(here $p_*$ is the norm on the dual space $E^*$, which is dual to $p$).
Examples of $\kappa$-regular norms with ``moderate'' $\kappa$
include the spaces $(\bR^n,\|\cdot\|_p)$ ($L_p$ on an $n$-point
set with unit masses of points) and the spaces $(\bR^{m\times
n},|\cdot|_p)$, $2\leq p\leq\infty$, of $m\times n$ matrices with
the Shatten norms $|X|_p=\|\sigma(X)\|_p$, $\sigma(X)$ being the
vector of singular values of a matrix $X$; in both cases,
$p\in[2,\infty]$. The spaces of the first series are
$\kappa$-regular with $\kappa=O(1)\min[p,\ln(n+1)]$, while the
spaces of the second series are $\kappa$-regular with
$\kappa=O(1)\min[p,\ln(m+1),\ln(n+1)]$.
\par
Norms $p(\cdot)$ satisfying (\ref{eq555}) play important role in
the theory of Banach spaces (where they are called norms with
smoothness modulus of power 2). In particular, a number of results on the properties of
martingales taking values in Banach spaces with smooth norms (see,
e.g., \cite{Garling1,Garling2}) are available. However, we were unable to
locate in the literature a result equivalent to Theorem \ref{mainthe}
which establishes  the validity of (somehow refined) bound $(*)$
in the case of a $\te^2$-regular space $(E,\|\cdot\|)$. Thus,
the main result of this paper, to the best of our (perhaps
incomplete)  knowledge, is new. The preliminary and slightly less
accurate, version of Theorem \ref{mainthe} was announced in
\cite{Nem2003} and proved in the preprint \cite{Nem2004OnLine}.
\par
While the question we
address seems to be important by its own right, our interest in it
stems mainly from various applications of (somehow rudimentary)
bounds of type $(*)$ we have encountered over the years. These applications include
investigating performance of Euclidean and non-Euclidean
stochastic approximation \cite{NYu,JudNem}, nonparametric
statistics \cite{NemPolTsy,JudNem,NemDiffIn}, optimization under
uncertainty \cite{Nem2003}, investigating quality of semidefinite
relaxations of some difficult combinatorial problems
\cite{Nem2007}, etc.
\par
Our paper is organized as follows:
the main result on large
deviations (Theorem \ref{mainthe}) is formulated in Section
\ref{mainresult}. Section \ref{sect1} contains instructive
examples and characterizations of $\kappa$-regular spaces, along
with a kind of ``calculus'' of these spaces. All proofs are placed
in the appendix.
\par
In what follows, if not explicitly stated otherwise, we suppose all the relations between random variables to hold a.s..

\section{Main result}
\label{mainresult}
\subsection{Regular spaces}
We start with the following
\begin{definition}\label{wwdef:smooth}
Let $(E,\|\cdot\|)$ be a finite-dimensional normed space and let
$\kappa\geq 1$.
\\
{\rm (i)} The
function $p(x)=\|x\|^2$ called $\kappa$-smooth if it is continuously differentiable and
\begin{equation}\label{wweq1}
\forall x,y\in E: p(x+y)\leq p(x) +Dp(x)[y] + \kappa p(y).
\end{equation}
\\
{\rm (ii)} Space $(E,\|\cdot\|)$ (and the norm $\|\cdot\|$ on $E$)
is called $\kappa$-regular, if there exists
$\kappa_+\in[1,\kappa]$ and a norm $\|\cdot\|_+$ on $E$ such that
$(E,\|\cdot\|_+)$ is $\kappa_+$-smooth and $\|\cdot\|_+$ is
$\kappa/\kappa_+$-compatible with $\|\cdot\|$, that is,
\begin{equation}\label{wweq2}
\forall x\in E: \|x\|^2\leq \|x\|_+^2\leq
{\kappa\over\kappa_+}\|x\|^2.
\end{equation}
{\rm (iii)} The constant $\kappa(E,\|\cdot\|)$ of regularity of
$E,\|\cdot\|$ is the infinum (clearly achievable) of those
$\kappa\geq1$ for which $(E,\|\cdot\|)$ is $\kappa$-regular.
\end{definition}
As an immediate example, an Euclidean space $(\bR^n,\|\cdot\|_2)$
is 1-smooth and thus 1-regular.
\subsection{Main result}
Assume that we are given\begin{itemize}
\item a finite-dimensional space $(E,\|\cdot\|)$, \item a Polish space $\Omega$ with Borel
  probability measure $\mu$,
   and
 \item
 a sequence ${\cF}_0=\{\emptyset,\Omega\}\subset{\cF}_1\subset{\cF}_2\subset...$ of $\sigma$-sub-algebras of the Borel
  $\sigma$-algebra of $\Omega$. \end{itemize}
  We denote by $\bE_i$, $i=1,2,...$
  the conditional expectation w.r.t. ${\cF}_i$, and by $\bE\equiv \bE_0$
  the expectation w.r.t. $\mu$.
  \par We further assume that we are given
an $E$-valued martingale-difference sequence
$\xi^\infty=\{\xi_i\}_{i=1}^\infty$ of Borel $E$-valued functions
on $\Omega$ such that $\xi_i$ is $\cF_i$-measurable and
$$
\bE_{i-1}\left\{\xi_{i}\right\}\equiv 0,\,i=1,2,...
$$
 An immediate consequence of Definition \ref{wwdef:smooth} of the regular norm is as follows:
assume that an $E$-valued martingale-difference
$\xi=\{\xi_t\}_{t=1}^\infty$ is square-integrable:
$\bE\left\{\|\xi_t\|^2\right\}\leq\sigma_t^2<\infty.$
Then
\[
\bE\left\{\|S_n\|^2\right\}\leq\kappa
\sum\limits_{t=1}^n\sigma_t^2.
\]
Indeed, $\|\cdot\|_+$ is $\kappa_+$-smooth, we have
$$
p(S_{t+1})\leq p(S_t)+Dp(S_t)[\xi_{t+1}]+\kappa_+ p(\xi_{t+1})
$$
whence, taking expectations and making use of the fact that $\xi$
is a martingale-difference,
$$
\bE\left\{p(S_{t+1})\right\}\leq\bE\left\{p(S_{t})\right\} +
\kappa_+
\bE\left\{p(\xi_{t+1})\right\}\leq\bE\left\{p(S_{t})\right\}+\kappa\bE\left\{\|\xi_{t+1}\|^2\right\}
$$
by the right inequality of (\ref{wweq2}). Then, by the left inequality of (\ref{wweq2}),
$$
\bE\left\{\|S_n\|^2\right\}\le \bE\left\{\|S_n\|_+^2\right\}\leq
\kappa\sum\limits_{t=1}^n\bE\left\{\|\xi_{t}\|^2\right\}\leq
\kappa\sum\limits_{t=1}^n\sigma_t^2.
$$
Our primary objective is to establish exponential bounds on the probabilities
of large deviations for an $E$-valued martingale difference
$\{\xi_i\}$. To this end, we impose on $\{\xi_i\}$ a ``light
tail'' assumption as follows.
Let $\alpha\in[1,2]$ and a sequence
$\sigma^\infty=\{\sigma_i>0\}_{i=1}^\infty$ of (deterministic)
positive reals be given. We introduce the following condition on the sequence $\xi^\infty$:
$$ \forall i\geq1:
\bE_{i-1}\left\{\exp\{\|\xi_i\|^\alpha\sigma_i^{-\alpha}\}\right\}\leq
\exp\{1\} \hbox{\ almost surely}
\eqno{(\cC_\alpha[\sigma^\infty])}
$$
Our main result is the large deviation bound for
$S_N=\sum_{i=1}^N\xi_i$ as follows:
\begin{theorem}\label{mainthe}
Let $(E,\|\cdot\|)$ be $\kappa$-regular, let $E$-valued
martingale-difference $\xi^\infty$ satisfy
$(\cC_\alpha[\sigma^\infty])$, and let $S_N=\sum_{i=1}^N \xi_i$,
$\sigma^N=[\sigma_1;...;\sigma_N]$. Then
\par {\rm (i)} for $1\leq\alpha\leq2$, one has for all $N\ge 1$ and $\g\ge 0$:
\begin{equation}\label{nfineq3}
\Prob\left\{\|S_N\|\geq
\left[\sqrt{2e\kappa}+\sqrt{2}\gamma\right]\sqrt{\sum_{i=1}^N\sigma_i^2}\right\}\leq
2\exp\left\{-{1\over
64}\min\left[\gamma^2;\gamma_*^{2-\alpha}\gamma^\alpha\right]\right\},
\end{equation}
where
\begin{equation}\label{Q.E.D.1.ini}
\gamma_*\equiv\gamma_*(\alpha,\nu^N)=\left\{\begin{array}{ll}
\begin{array}{l}32\left[{8\alpha_*\over
2^{\alpha_*}}\right]^{{\alpha-1\over2-\alpha}}\left[{\|\nu^N\|_2\over\|\nu^N\|_{\alpha_*}}
\right]^{{\alpha\over2-\alpha}} \geq
16\left[{\|\nu^N\|_2\over\|\nu^N\|_{\alpha_*}}
\right]^{{\alpha\over2-\alpha}}\geq16,\\
\left[\alpha_*={\alpha\over\alpha-1},\nu^N=[\nu_1;...;\nu_N]\right],\\
\end{array}&1<\alpha<2,\\
\lim_{\alpha \to
1+0}\gamma_*(\alpha,\nu^N)=16{\|\nu^N\|_2\over\|\nu^N\|_\infty},&\alpha=1,\\
\lim_{\alpha\to 2-0}\gamma_*(\alpha,\nu^N)=+\infty,&\alpha=2.\\
\end{array}\right.
\end{equation}
\par {\rm (ii)} When $\alpha=2$, the bound {\rm (\ref{fineq3})}
improves to
\begin{equation}\label{nfineq3_impr}
(\forall N\geq1,\gamma\geq0): \Prob\left\{\|S_N\|\geq
\left[\sqrt{2\kappa}+\sqrt{2}\gamma\right]\sqrt{\sum_{i=1}^N\sigma_i^2}\right\}\leq\exp\{-\gamma^2/3\}.
\end{equation}
\par
{\rm (iii)} When  the condition
$\bE_{i-1}\left\{\exp\{\|\xi_i\|^2\sigma_i^{-2}\}\right\}\leq\exp\{1\}$
in $(\cC_2[\sigma^\infty])$ is strengthened to
$\|\xi_i\|\leq\sigma_i$ almost surely, $i=1,2,...,$ the bound {\rm
(\ref{fineq3})} improves to
\begin{equation}\label{nfineq3_impr_impr}
(\forall N\geq1,\gamma\geq0): \Prob\left\{\|S_N\|\geq
\left[\sqrt{2\kappa}+\sqrt{2}\gamma\right]\sqrt{\sum_{i=1}^N\sigma_i^2}\right\}\leq
\exp\left\{-\gamma^2/2\right\}.
\end{equation}
\end{theorem}
\section{Regular spaces}
To make Theorem \ref{mainthe} meaningful, we need to point out a
spectrum of interesting $\kappa$-smooth/regular spaces, and this
is the issue we consider in this Section.
\subsection{Basic examples}
\label{sect1} Let $E$ be an $n$-dimensional linear space, and let
$\|\cdot\|$ be a norm on $E$. It is well known \cite{FritzJohn}
that there exists an ellipsoid $Q$ centered at the origin such
that $Q\subset \{x\in E:\|x\|\leq1\}\subset\sqrt{n}Q$, or,
equivalently, there exists a Euclidean norm $\|\cdot\|_+$ on $E$
such that $\|x\|^2\leq\|x\|_+^2\leq n\|x\|^2$. Since the Euclidean
space $(E,\|\cdot\|_+)$ is $1$-smooth, we conclude that
\begin{example}\label{wwex:R_n}.  Every
finite-dimensional  normed space $(E,\|\cdot\|)$ is $(\dim
E)$-regular.\end{example}\par We are about to present a number of
less trivial examples, those where  the regularity parameter
$\kappa$ is dimension-independent (or nearly so).
\begin{example}\label{wwex:l_p} Let $2\leq p\leq\infty$. The space
$(\bR^n,\|\cdot\|_p)$ with $n\geq3$ is $\kappa_p(n)$-regular with
\begin{equation}\label{wwkappa_p(n)}
\kappa_p(n)=\min\limits_{{2\leq\rho\leq p\atop \rho<\infty}}
(\rho-1)n^{{2\over\rho}-{2\over p}}\leq \min[p-1,2\ln(n)]
\end{equation}
\end{example}
\begin{example}\label{wwex:matrices} Let $2\leq p\leq\infty$.
The norm $|X|_p=\|\sigma(X)\|_p$ on the space $\bR^{m\times n}$ of
$m\times n$ real matrices, where $\sigma(X)$ is the vector of
singular values of $X$, is $\kappa_p(m,n)$-regular, with
\begin{equation}\label{wwkappa_p(m,n)}
\begin{array}{l}
\kappa_p(m,n)= \min\limits_{{2\leq\rho<\infty\atop \rho\leq p}}
\max[2,\rho-1](\min(m,n))^{{2\over\rho}-{2\over p}}\\
\leq \min\left[\max[2,p-1], (2\ln (\min[m,n]+2)
-1)\exp\{1\}\right].\\
\end{array}
\end{equation} \end{example}
The proof of the bound (\ref{wwkappa_p(m,n)}) is based upon the fact which is
important by its own right:
\begin{proposition}\label{prop2.1}
Let $\Delta$ be an open interval on the axis, and $f$ be a C$^2$
function on $\Delta$ such that for certain
$\theta_\pm,\mu_\pm\in\bR$ one has
\begin{equation}\label{eqno12}
\forall (a<b,a,b\in\Delta): \theta_-{f''(a)+f''(b)\over 2}+\mu_-
\leq {f'(b)-f'(a)\over b-a}\leq\theta_+{f''(a)+f''(b)\over
2}+\mu_+
\end{equation}
Let, further, ${\cal  X}_n(\Delta)$ be the set of all $n\times n$
symmetric matrices with eigenvalues belonging to $\Delta$. Then
${\cal X}_n(\Delta)$ is an open convex set in the space $\bS^n$ of
$n\times n$ symmetric matrices, the function
$$ F(X) = \Tr(f(X)) : {\cal X}_n(\Delta) \to\bR
$$ is C$^2$, and for every $X\in{\cal X}_n(\Delta)$ and every
$H\in\bS^n$ one has
\begin{equation}\label{eqno13}
\theta_-\Tr(Hf''(X)H)+\mu_-\Tr(H^2)\leq D^2F(X)[H,H]\leq
\theta_+\Tr(Hf''(X)H)+\mu_+\Tr(H^2).
\end{equation}
\end{proposition}
\subsection{Dual characterization of smoothness and regularity}
The following well-known fact can be seen as  dual
characterization of $\kappa$-smoothness:
\begin{proposition} \label{char} Let
$(E,\|\cdot\|)$ be a finite-dimensional normed space, $E^*$ be the
space dual to $E$, $\|\cdot\|_*$ be the norm on $E^*$ dual to
$\|\cdot\|$; and let $\langle\xi,x\rangle$ stand for the value of a linear form
$\xi\in E^*$ on a vector $x\in E$. Let also $f(x)={1\over
2}\|x\|^2:E\to\bR$ and $f_*(\xi)={1\over 2}\|\xi\|_*^2:E^*\to
\bR$. The following 
properties are equivalent to each other:
\begin{itemize} \item[{\rm (i)}] $(E,\|\cdot\|)$ is
$\kappa$-smooth; \item[{\rm (ii)}] $\partial f(x)=\{f'(x)\}$ is a
singleton for every $x$, and
\begin{equation}\label{beq1}
\langle f'(x)-f'(y),x-y\rangle \leq \kappa\|x-y\|^2\;\;\forall
x,y\in E;
\end{equation}
\item[{\rm (iii)}] $f$ is continuously differentiable, and
$f'(\cdot)$ is Lipschitz continuous with constant $\kappa$:
\begin{equation}\label{beq2}
\|f'(x)-f'(y)\|_*\leq\kappa\|x-y\|\;\;\forall x,y\in E;
\end{equation}
\item[{\rm (iv)}] One has
$$
\forall (\xi,\eta\in E^*,\;x\in\partial f_*(\xi),\;y\in\partial
f_*(\eta)):\;\; \langle\xi-\eta,x-y\rangle \geq
\kappa^{-1}\|\xi-\eta\|_*^2;\\
$$
\item[{\rm (v)}] One has
$$
\forall (\xi,\eta\in E^*,\;x\in \partial f_*(\xi),\;y\in\partial
f_*(\eta)): \;\;\|x-y\|\geq \kappa^{-1}\|\xi-\eta\|_*;
$$
\item[{\rm (vi)}] One has
$$
\forall (\xi,\eta\in E_*,\;x\in\partial f_*(\xi)):\;\; f_*(\xi+\eta)\geq
f_*(\xi)+\langle\eta,x\rangle +{1\over 2\kappa}\|\eta\|_*^2.
$$
\end{itemize}
\end{proposition}
Another characterization of regular spaces is as follows:
\begin{proposition}\label{propagain}\label{prop27} Let
$(E,\|\cdot\|)$ be a finite-dimensional normed space, $E^*$ be the space
dual to $E$, $\|\cdot\|_*$ be the norm on $E^*$ dual to
$\|\cdot\|$, and  let $\langle\xi,x\rangle$ stand for the value of a linear
form $\xi\in E^*$ on a vector $x\in E$. Let also $B_*$ be the
unit $\|\cdot\|_*$-ball of $E^*$.
\par
{\rm (i)} If $(E,\|\cdot\|)$ is $\kappa$-regular, then the exists
a continuous function $V:B_*\to\bR$ which is strongly convex, with
coefficient $1$ w.r.t. $\|\cdot\|_*$, on $B_*$, that is, possesses
the following equivalent to each other properties:
\begin{equation}\label{equivalent}
\begin{array}{ll}(a)&\forall (\xi,\;\eta\in\inter B_*,\;x\in\partial v(\xi),y\in\partial v(\eta)): \;\;\langle
\xi-\eta,x-y\rangle \geq \|\xi-\eta\|_*^2,\\
(b)&\forall (\xi,\eta:\;\xi,\xi+\eta\in\inter B_*,\;x\in\partial
v(\xi)): \;\;v(\xi+\eta)\geq v(\xi)+\langle \eta,x\rangle+{1\over
2}\|\eta\|_*^2;\\
\end{array}
\end{equation}
and, in addition, is such that
\begin{equation}\label{issuchthat}
\max\limits_{B_*} v-\min\limits_{B_*} v\leq {\kappa\over2}
\end{equation}
\par
{\rm (ii)} Assume that the unit ball $B_*$ of $(E^*,\|\cdot\|_*)$
admits a function $v$ satisfying {\rm (\ref{equivalent}),
(\ref{issuchthat})}. Then $(E,\|\cdot\|)$ is $O(1)\kappa$-regular
with an appropriately chosen absolute constant $O(1)$.
\end{proposition}
\subsection{``Calculus'' of smooth and regular spaces}
\begin{proposition}\label{prop28}
Let $(E,\|\cdot\|_E)$ be a finite-dimensional normed space, $L$ be
a linear subspace of $E$, and $F=E/L$ be the factor-space of $E$
equipped with the factor-norm $\|\bar{f}\|_F=\min_{f\in
\bar{f}}\|f\|_E$. If $(E,\|\cdot\|_E)$ is $\kappa$-smooth
($\kappa$-regular), then $(L,\|\cdot\|_E)$ and $(F,\|\cdot\|_F)$
also are $\kappa$-smooth, respectively, $\kappa$-regular.
\end{proposition}
\begin{proposition} \label{calculus}  {\rm
(i)} Let $p\in[2,\infty]$, and let $(E_i,\|\cdot\|_i)$ be
finite-dimensional $\kappa$-smooth spaces, $i=1,...,m>2$. The
space $E=E_1\times ...\times E_m$ equipped with the norm
$$
\|(x^1,...,x^m)\|=\left(\sum\limits_{i=1}^m
\|x^i\|_i^p\right)^{1/p}
$$
(the right hand side is $\max\limits_i\|x^i\|_i$ when $p=\infty$)
is $\kappa^+$-regular with
\begin{equation}\label{kappa_plus}
\kappa^+=\min\limits_{2\leq \rho\leq p}
[\kappa+\rho-1]m^{{2\over\rho}-{2\over
p}}\leq\min[\kappa+p-1,[\kappa+2\ln(m)-1]\exp\{1\}].
\end{equation}
\par
{\rm (ii)} Let $\|\cdot\|_i$ be $\kappa$-smooth norms on $E$. Then
the norm
$$
\|x\|=\sum\limits_{i=1}^m \|x\|_i
$$
is $m\kappa$-regular on $E$.
\end{proposition}
\begin{proposition}\label{calculusI}
 {\rm (i)} Let $p\in[2,\infty]$, and let $(E_i,\|\cdot\|_i)$ be
finite-dimensional $\kappa$-regular spaces, $i=1,...,m>2$. The
space $E=E_1\times ...\times E_m$ equipped with the norm
$$
\|(x^1,...,x^m)\|=\left(\sum\limits_{i=1}^m
\|x^i\|_i^p\right)^{1/p}
$$
(the right hand side is $\max\limits_i\|x^i\|_i$ when $p=\infty$)
is $\kappa^{++}$-regular with
\begin{equation}\label{kappa_plus_plus}
\kappa^{++}=2\min\limits_{2\leq \rho\leq p}
[\kappa+\rho-1]m^{{2\over\rho}-{2\over
p}}\leq2\min[\kappa+p-1,[\kappa+2\ln(m)-1]\exp\{1\}].
\end{equation}
\par
{\rm (ii)} Let $\|\cdot\|_i$ be $\kappa$-regular norms on a
finite-dimensional space $E$. Then the norm
$$
\|x\|=\sum\limits_{i=1}^m \|x\|_i
$$
is $2m\kappa$-regular on $E$.
\end{proposition}
\section{Appendix: Proofs}
\subsection{Proofs for Section \ref{sect1}}
\subsubsection{Justifying the Examples}
\paragraph{Example \ref{wwex:l_p}:} Let $2\leq \rho<\infty$. We claim that in this case the space
$(R^n,\|\cdot\|_\rho)$ is $(\rho-1)$-smooth. Indeed, the function
$p(x)=\|\cdot\|_\rho^2$ is convex, continuously differentiable
everywhere and twice continuously differentiable outside of the
origin; for such a function, (\ref{wweq1}) holds true if and only
if
\begin{equation}\label{wweq3}
D^2p(x)[h,h]\leq 2\kappa_+ p(h)\,\,\forall (x,h\in E, x\neq0);
\end{equation}
since $p(\cdot)$ is homogeneous of degree 2, the validity of
(\ref{wweq3}) for all $x,h$ is equivalent to the validity of the
relation for all $h$ and all $x$ normalized by the requirement
$p(x)=1$. Given such an $x$ and $h$ and assuming $\rho>2$, we have
$$
\begin{array}{rcl}
Dp(x)[h]&=&2\left(\sum\limits_i|x_i|^\rho\right)^{{2\over
\rho}-1}\sum\limits_i |x_i|^{\rho-1}\sign(x_i)h_i\\
D^2p(x)[h,h]&=&2\underbrace{\left({2\over\rho}-1\right)}_{\leq0}\left(\sum\limits_i|x_i|^\rho\right)^{{2\over
\rho}-2}\left(\sum\limits_i |x_i|^{\rho-1}\sign(x_i)h_i\right)^2\\
&&+2\big(\underbrace{\sum\limits_i|x_i|^\rho}_{=1}\big)^{{2\over
\rho}-1} \sum\limits_i(\rho-1)|x_i|^{\rho-2}h_i^2
\leq 2(\rho-1)\sum\limits_i|x_i|^{\rho-2}h_i^2\\
&\leq&
2(\rho-1)\left(\sum\limits_i(|x_i|^{\rho-2})^{{\rho\over\rho-2}}\right)^{{\rho-2\over\rho}}\left(\sum\limits_i(|h_i|^2)^{{\rho\over
2}}\right)^{{2\over\rho}}\\
&=&2(\rho-1)\|h\|_\rho^2=2(\rho-1)p(h)\\
\end{array}
$$
as required in (\ref{wweq3}) when $\kappa_+=\rho-1$. In the case
of $\rho=2$ relation (\ref{wweq3}) with $\kappa_+=\rho-1=1$ is
evident.
\par
Now, when $\rho\in[2,p]$ and $x\in\bR^n$, one has
$\|x\|^2_\rho/\|x\|_p^2\in[1,n^{{2\over\rho}-{2\over p}}]$, so
that $(\bR^n,\|\cdot\|_p)$ is $\kappa$-regular with
$\kappa=(\rho-1)n^{{2\over\rho}-{2\over p}}$, and
(\ref{wwkappa_p(n)}) follows. \blacksquare
\paragraph{Example \ref{wwex:matrices}:}
1$^0$. We start with the following
\begin{lemma}\label{wweven} Let
$\rho\geq2$. Then the space $\bS^n$ of symmetric $n\times n$
matrices with the norm $|X|_\rho$ is $\kappa$-smooth with
\begin{equation}\label{wwkappa} \kappa=\max[2,\rho-1]. \end{equation}
\end{lemma}
{\bf Proof.} The statement is evident when $\rho=2$; thus, from
now on we assume that $\rho>2$.  Let us apply Proposition
\ref{prop2.1} to $\Delta=\bR$, $f(t)=|t|^\rho$ with
$\theta_-=\mu_-=0$, $\mu_+=0$ and $\theta_+= \max\left[{2\over
\rho-1},1\right]$ (this choice, as it is immediately seen,
satisfies (\ref{eqno12})). By Proposition, the function $F(X) =
|X|_\rho^\rho$ on $\bS^n$ is twice continuously differentiable,
and
\begin{equation}\label{eqno14}\forall X,H: 0\leq
D^2F(X)[H,H]\leq\theta_+
\Tr(f''(x)H^2),\,\,\theta_+=\max\left[{2\over \rho-1},1\right].
\end{equation}
  It follows that the function $p(X) = |X|_\rho^2= (F(X))^{{2\over\rho}}$
is continuously differentiable everywhere and twice continuously
differentiable outside of the origin. For $X\neq0$ we have
$Dp(X)[H]={2\over\rho}(F(X))^{{2\over\rho}-1}DF(X)[H]$, whence
\begin{equation}\label{eqno15}
\begin{array}{l}
X\neq0\Rightarrow
D^2p(X)[H,H]={2\over\rho}\underbrace{\left[{2\over\rho}-1\right]}_{<0}
(F(X))^{{2\over\rho}-2}(DF(X)[H])^2+{2\over\rho}(F(X))^{{2\over\rho}-1}
D^2F(X)[H,H]\\
\leq{2\over\rho}(F(X))^{{2\over\rho}-1}\theta_+ \Tr(f''(x)H^2).
\\
\end{array}
\end{equation}
Setting $Z={1\over\rho(\rho-1)}(F(X))^{{2\over\rho}-1}f''(X)$,
$p={\rho\over\rho-2}$, it is immediately seen that $|Z|_p=1$. From
(\ref{eqno15}) we have
\begin{equation}\label{eqno16}
\begin{array}{l}
D^2p(X)[H,H]\leq 2\Theta_+(\rho-1)\Tr(ZH^2)\leq
2\theta_+(\rho-1)|Z|_p|H^2|_{{p\over p-1}}=
2\theta_+(\rho-1)|H^2|_{{\rho\over2}}\\
=2\theta_+(\rho-1)|H|_\rho^2.\\
\end{array}
\end{equation}
 Now, if $X, Y\in\bS^n$  are such that the segment $[X;X + Y ]$ does
not contain the origin, then
$$
\exists \gamma\in(0,1): p(X+Y)\leq p(X)+Dp(X)[Y]+{1\over
2}D^2p(X+\gamma Y)[Y,Y],
$$
and (\ref{eqno16}) implies that for the outlined $X,Y$ one has
$$
p(X+Y)\leq p(X)+Dp(X)[Y]+\theta_+(\rho-1)p(Y).
$$
Since $p$ is C$^1$, the resulting inequality, by continuity, is
valid for all $X,Y$. \blacksquare
\par
2$^0$. Now we can complete the justification of Example
\ref{wwex:matrices}. W.l.o.g. we may assume that $m\leq n$. Given
an $m\times n$ matrix $X$, let
$S(X)=\left[\begin{array}{c|c}&X\cr\hline
X^T&\cr\end{array}\right]\in \bS^{m+n}$. One clearly has
$$
\|\sigma(X)\|_\rho=|X|_\rho=2^{-1/\rho}|S(X)|_\rho,
$$
whence, by Lemma \ref{wweven} and due to the fact that the mapping
$X\mapsto S(X): \bR^{m\times n}\to \bS^{m+n}$ is linear, the norm
$|\cdot|_\rho$, treated as a norm on $\bR^{m\times n}$, is
$\max[2,\rho-1]$-smooth whenever $\rho\geq2$. Since
$\sigma(X)\in\bR^m$ for $X\in\bR^{m\times n}$, for every
$\rho\in[2,\infty)$ such that $\rho\leq p$ one has
$$
|X|_p^2\leq |X|_\rho^2\leq m^{{2\over\rho}-{2\over p}}|X|_p^2.
$$
Thus,  the space $(\bR^{m\times n},|\cdot|_p)$ is $\kappa$-regular
with $ \kappa =\min\limits_{{2\leq\rho<\infty\atop \rho\leq p}}
\max[2,\rho-1]m^{{2\over\rho}-{2\over p}},
$
and we arrive at (\ref{wwkappa_p(m,n)}). \blacksquare
\subsubsection{Proof of Proposition \ref{prop2.1}}
Let $\{f_k(t)\}$ be a sequence of polynomials converging to $f$,
along with the first and the second derivatives, uniformly on
every compact subset of $\Delta$. For a polynomial
$p(t)=\sum_{j=0}^N p_jt^j$ the function $P(X)=\Tr(\sum_jp_jX^j)$
is a polynomial on $\bS^n$. Let now $X,H\in \bS^n$, let
$\lambda_s=\lambda_s(X)$ be the eigenvalues of $X$,
$X=U\Diag\{\lambda\}U^T$ be the eigenvalue decomposition of $X$,
and let $\widehat{H}$ be such that $H=U\widehat{H}U^T$. We have
\begin{equation}\label{eqno89}
\begin{array}{rclr}
P(X)&=&\sum_{s=1}^n p(\lambda_s(X))&(a)\\
DP(X)[H]&=&\Tr(\sum_{j=1}^N\sum_{s=0}^{N-1}X^sHX^{N-s-1}=\Tr(p'(X)H)=\sum_{s=1}^n
p'(\lambda_s(X))\widehat{H}_{ss} &(b)\\
\end{array}
\end{equation}
Further, let $\gamma$ be a closed contour in the complex plane
encircling all the eigenvalues of $X$. Then
$$
\begin{array}{l}
DP(X)[H]=\Tr(p'(X)H)={1\over2\pi\imath}\oint\limits_\gamma
p'(z)\Tr((zI-X)^{-1}H)dz\\
\Rightarrow D^2P(X)[H,H]={1\over 2\pi\imath}\oint\limits_\gamma
p'(z)\Tr((zI-X)^{-1}H(zI-X)^{-1}H)dz= {1\over2\pi\imath}
\oint\limits_\gamma \sum_{s,t=1}^n
{\widehat{H}_{st}^2p'(z)\over(z-\lambda_s)(z-\lambda_t)}dz.
\end{array}
$$
Computing the residuals, we get
\begin{equation}\label{eqno90}
D^2P(X)[H,H]=\sum_{s,t}\Gamma_{s,t}[p]\widehat{H}_{st}^2,\quad
\Gamma_{s,t}[p]=\left\{\begin{array}{ll}
{p'(\lambda_s)-p'(\lambda_t)\over
\lambda_s-\lambda_t},&\lambda_s\neq\lambda_t\\
p''(\lambda_s),&\lambda_s=\lambda_t\\
\end{array}\right.
\end{equation}
Substituting $p=f_k$ into (\ref{eqno89}.$a,b$) and (\ref{eqno90}),
we see that the sequence of polynomials $F_k(X)=\Tr(f_k(X))$
converges, along with the first and the second order derivatives,
uniformly on compact subsets of ${\cal X}_n(\Delta)$; by
(\ref{eqno89}.$a$), the limiting function is exactly $F(X)$. We
conclude that $F(X)$ is C$^2$ on ${\cal X}_n(\Delta)$ and that the
first and the second derivatives of this function are limits, as
$k\to\infty$, of the corresponding derivatives of $F_k(X)$, so
that for $X=U\Diag\{\lambda\}U^T\in{\cal X}_n(\Delta)$ (where $U$
is orthogonal) and every $H=U\widehat{H}U^T\in\bS^n$ we have
\begin{equation}\label{eqno91}
\begin{array}{rcl}
DF(X)[H]&=&\sum_s f'(\lambda_s) \widehat{H}_{ss}=\Tr(f'(X)H)\\
D^2F(X)[H,H]&=&\sum_{s,t}\Gamma_{s,t}[f]\widehat{H}^2_{st}\\
\end{array}\end{equation}
So far, we did not use (\ref{eqno12}). Invoking the right
inequality in (\ref{eqno12}), we get
$$
\begin{array}{l}
D^2F(X)[H,H]\leq\sum_{s,t}\left[\theta_+{f''(\lambda_s)+
f''(\lambda_t)\over2} +\mu_+\right]\widehat{H}^2_{st}=
\theta_+\sum_sf''(\lambda_s)\sum_t\widehat{H}^2_{st}+\mu_+
\sum_{s,t}\widehat{H}^2_{st}\\
=\theta_+\Tr(\Diag\{f''(\lambda_1),...,f''(\lambda_n)\}\widehat{H}^2)
+\mu_+\Tr(\widehat{H}^2)=\theta_+\Tr(f''(X)H^2)+\mu_+\Tr(H^2),
\\
\end{array}
$$
which is the right inequality in (\ref{eqno13}). The derivation of
the left inequality in (\ref{eqno13}) is similar. \blacksquare
\subsubsection{Proof of Proposition \ref{char}}
\paragraph{\bf (i)$\Rightarrow$(iii)}: We are in the situation when $f$
is continuously differentiable. Convolving $f(\cdot)$ with smooth
nonnegative kernels $\delta_k(\cdot)$
 with unit integral and support shrinking to
origin as $k\to\infty$, we get a sequence $f_k(\cdot)$ of smooth
functions converging to $f(\cdot)$, along with first order
derivatives, uniformly on compact sets. We have
$$
\begin{array}{rcl}
f_k(x+y)&=&\int f(x-z+y)\delta(z)dz \leq \int
[f(x-z)+\langle f'(x-z),y\rangle +\kappa f(y)]\delta(z)dz\\
&=&f_k(x)+\langle f_k^\prime(x),y\rangle + \kappa f(y)\\
\end{array}
$$
From the resulting inequality combined with smoothness and
convexity of $f_k$ it follows that
$$
0\leq D^2f_k(x)[h,h]\leq \kappa \|h\|^2\,\,\forall x,h\in E.
$$
Thus, if $\|h\|=\|d\|=1$, then
$$4D^2f_k(x)[h,d]=D^2f_k(x)[h+d,h+d]-D^2f_k(x)[h-d,h-d]\leq
\kappa\|h+d\|^2\leq 4\kappa$$. Whence $D^2f_k(x)[h,d]\leq \kappa$
whenever $\|h\|=\|d\|=1$, or, which is the same by homogeneity,
$$
|D^2f_k(x)[h,d]|\leq \kappa \|h\|\|d\| \,\,\forall x,h,d.
$$
Consequently,
$$
|\langle f_k^\prime(y)-f_k^\prime(x),h\rangle|=|\int\limits_0^1
D^2f_k(x+t(y-x))[y-x,h]dt|\leq \int\limits_0^1
\kappa\|y-x\|\|h\|dt \leq \kappa \|y-x\|\|h\|,
$$
whence, taking maximum over $h$ with $\|h\|=1$,  $$\|
f_k^\prime(y)-f_k^\prime(x)\|_*\leq \kappa\|y-x\|$$. As
$k\to\infty$, $f_k^\prime(x)$ converge to $f'(x)$, and we conclude
that $f'(\cdot)$ possesses the required Lipschitz continuity.
\myqed
\paragraph{(iii)$\Rightarrow$(ii):} evident
\paragraph{(ii)$\Rightarrow$(i):} A convex function on $\bR^n$ with a singleton differential at every point
 clearly is
 continuously differentiable, so that in the case of (ii) $f$ is continuously differentiable. Besides this, in the case of (ii)
 we have
$$
\begin{array}{rcl}
f(x+y)&=&f(x)+\langle f'(x),y\rangle + \int\limits_0^1 \langle
f'(x+ty)-f'(x),y\rangle dt\\
&\leq&f(x)+\langle f'(x),y\rangle + \int\limits_0^1\kappa
t\|y\|^2dt=f(x)+\langle f'(x),y\rangle + \kappa f(y),\\
\end{array}
$$
which immediately implies (\ref{wweq1}) (recall that
$\|\cdot\|^2=2f(\cdot)$). \myqed
\paragraph{(iii)$\Leftrightarrow$(v):} The functions
$f(\cdot)$, $f_*(\cdot)$ are the Legendre transforms of each
other, so that $x\in\partial f_*(\xi)$ if and only if
$\xi\in\partial f(x)$. Now let (iii) be the case, and let
$\xi,\eta\in E^*$ and $x\in\partial f_*(\xi)$, $y\in\partial
f_*(\eta)$. Then $\xi=f'(x)$, $\eta=f'(y)$ and therefore, due to
(iii), $$ \|\xi-\eta\|_*\leq \kappa\|x-y\|,$$ so that (v) takes
place. Vice versa, let (v) take place, and let $x,y\in E$,
$\xi\in\partial f(x)$, $\eta\in\partial f(y)$. Then $x\in\partial
f_*(\xi)$, $y\in\partial f_*(y)$, and therefore (v) says that
$$\|\xi-\eta\|_*\leq\kappa\|x-y\|.$$ We conclude that if $x=y$, then
$\xi=\eta$, that is, $\partial f(x)$ always is a singleton,
meaning that $f$ is continuously differentiable, and that the
inequality in (iii) takes place, that is, (iii) holds true. \myqed
\paragraph{(iv)$\Leftrightarrow$(iii):} Let (iv) take place. If there exists $x\in E$ such that $\partial f(x)$ is not a singleton, then, choosing $\xi,\eta\in\partial f(x)$
with $\xi\neq\eta$, we would have $x\in\partial f_*(\xi)$,
$x\in\partial f_*(\eta)$, whence by (iv) we should have $$\langle
\xi-\eta,x-x\rangle \geq \kappa^{-1}\|\xi-\eta\|_*^2,$$ which is
impossible. Thus, $\partial f(x)$ is a singleton for every $x$, so
that $f$ is continuously differentiable. Besides this, with
$x,y\in E$ and $\xi=f'(x)$, $\eta=f'(y)$ we have $x\in\partial
f_*(\xi)$, $y\in\partial f_*(\eta)$, whence, by (iv), $$\langle
\xi-\eta,x-y\rangle\geq\kappa^{-1}\|\xi-\eta\|_*^2.$$ Since
$$\langle \xi-\eta,x-y\rangle\leq\|\xi-\eta\|_*\|x-y\|,$$ we get
$$\|\xi-\eta\|_*\|x-y\|\geq \kappa^{-1}\|\xi-\eta\|_*^2,$$ whence
$$\|\xi-\eta\|_*=\|f'(x)-f'(y)\|_*\leq\kappa\|x-y\|,$$ and thus
(iii)  takes place.
\par
Now let (iii) take place, and let us prove that (iv) takes place
as well, or, which is the same in the case of (iii), that $\langle
f'(x)-f'(y),x-y\rangle\geq\kappa^{-1}\|f'(x)-f'(y)\|^2$. Setting
$$g(u)=f(u)-\langle f'(y),u-y\rangle,$$ we get a continuously
differentiable convex function on $E$ such that
$$\|g'(x)-g'(y)\|_*\leq\kappa\|x-y\|$$ and $g'(y)=0$. Due to these
relations, $$g(y+h)\leq g(y)+{\kappa\over 2}\|h\|^2$$ for all $h$.
Now let $e\in E$ be such that $\langle g'(x),e\rangle=\|g'(x)\|_*$
and $\|e\|=1$. Due to $$\|g'(u)-g'(v)\|_*\leq \kappa\|u-v\|,$$ we
have
\bse
g(x-{\|g'(x)\|_*\over \kappa}e)&\leq& g(x)-\langle
g'(x),{\|g'(x)\|_*\over \kappa}e\rangle +{\kappa\over
2}\|{\|g'(x)\|\over\kappa}e\|^2\\&=&g(x)-{\|g'(x)\|_*^2\over
\kappa}+{\|g'(x)\|_*^2\over 2\kappa}=g(x)-{\|g'(x)\|_*^2\over
2\kappa}.
\ese
 On the other hand, $g$ attains its global minimum at
$y$, so that $$g(x)-{\|g'(x)\|_*^2\over 2\kappa}\geq
g(x-{\|g'(x)\|_*\over \kappa}e)\geq g(y).$$ We now have
\bse
g(y)+{\kappa\over 2}\|h\|^2&\geq& g(y+h)\geq g(x)+\langle
g'(x),y+h-x\rangle\\&\geq& g(y)+{\|g'(x)\|_*^2\over2\kappa}+\langle
g'(x),y+h-x\rangle,
\ese
whence
$$
\langle g'(x),x-y\rangle\geq {\|g'(x)\|_*^2\over2\kappa}+\langle
g'(x),h\rangle - {\kappa\over 2}\|h\|^2.
$$
This inequality is valid for all $h$; setting
$h={\|g'(x)\|_*\over\kappa}e$, the right hand side becomes
${\|g'(x)\|_*^2\over\kappa}$. Thus, $$\langle
f'(x)-f'(y),x-y\rangle =\langle g'(x),x-y\rangle\geq
{\|g'(x)\|_*^2\over\kappa}={\|f'(x)-f'(y)\|_*^2\over \kappa}.\;\;\;\hfill{\mbox{\myqed}}$$
\paragraph{(iv)
$\Rightarrow$(vi):} Let (iv) take place, let $\xi,\eta\in E^*$ and
$x\in \partial f_*(\xi)$. Setting $\xi_t=\xi+t\eta$,
$\phi(t)=f_*(\xi_t)$, $0\leq t\leq 1$, we get an absolutely
continuous function on $[0,1]$ with the derivative which is almost
everywhere given by $\phi'(t)=\langle \eta,x_t\rangle$, with
$x_t\in\partial f_*(\xi_t)$. We have
\bse
f_*(\xi+\eta)&=&\phi(1)=\phi(0)+\int_0^1\phi'(t)dt\\
&=&
\phi(0)+\int_0^1\langle \eta,x_t\rangle dt=\phi(0)+\int_0^1
[\langle\eta,x\rangle+\langle
\eta,x_t-x\rangle]dt\\
&=&\phi(0)+\langle\eta,x\rangle +\int_0^1t^{-1}\langle
(\xi+t\eta)-\xi,x_t-x\rangle dt \\&\geq& \phi(0)+\langle\eta,x\rangle
+\int_0^1t^{-1}\kappa^{-1}\|[\xi+t\eta]-\xi\|_*^2dt\\
&=&\phi(0)+\langle\eta,x\rangle +{1\over
2\kappa}\|\eta\|_*^2=f_*(\xi)+\langle\eta,x\rangle
+{1\over 2\kappa}\|\eta\|_*^,\ese
where the inequality is given by (iv). We end up with the
inequality required in (vi). \myqed
\paragraph{(vi)$\Rightarrow$(i):} Let (vi) be the case,
let $x\in E$ and $\xi\in \partial f(x)$, so that $x\in\partial
f_*(\xi)$. We have
\bse
f(x+y)&=&\max_{\eta\in E^*}\left[\langle \xi+\eta,x+y\rangle
-f_*(\xi+\eta)\right]\\
&\leq&\max_{\eta\in
E^*}\left[\langle\xi+\eta,x+y\rangle
-f_*(\xi)-\langle\eta,x\rangle -{1\over
2\kappa}\|\eta\|_*^2\right]\\
&=&\max_{\eta\in E^*}\left[\langle \xi,x+y\rangle
+\langle\eta,y\rangle - f_*(\xi)-{1\over
2\kappa}\|\eta\|_*^2\right]\\
&=&
\underbrace{\langle\xi,x\rangle-f_*(\xi)}_{f(x)}+\langle\xi,y\rangle
+\max_{\eta}\left[\langle\eta,y\rangle-{1\over2\kappa}\|\eta\|_*^2\right]
=f(x)+\langle\xi,y\rangle+{\kappa\over 2} \|y\|^2.
\ese
This relation along with the relation $f(x+y)\geq
f(x)+\langle\xi,y\rangle$ implies that $\xi$ is the Frechet
derivative of $f$ at $x$, whence $f$ is convex and differentiable,
and thus -- continuously differentiable function on $E$ which
satisfies the inequality
$$
f(x+y)\leq f(x)+\langle f'(x),y\rangle +{\kappa\over
2}\|y\|^2.\eqno{\hbox{\myqed}}
$$
\par We have proved that
(i)$\Leftrightarrow$(ii)$\Leftrightarrow$(iii)$\Leftrightarrow$(iv)$\Leftrightarrow$(v)
and (iv)$\Rightarrow$(vi)$\Rightarrow$(i), meaning that all 6
properties in question are equivalent to each other. \blacksquare
\subsubsection{Proof of Proposition \ref{propagain}}
\paragraph{(i):} Let $(E,\|\cdot\|)$ be $\kappa$-regular, and let
$\kappa_+\in[1,\kappa]$ and $\|\cdot\|_+$ be such that
$(E,\|\cdot\|_+)$ is $\kappa$-smooth and (\ref{wweq2}) holds true,
and let $\|\cdot\|_{+,*}$ be the norm on $E^*$ dual to
$\|\cdot\|_+$; note that
\begin{equation}\label{notethat}
{\kappa_+\over\kappa}\|\cdot\|_*^2\leq \|\cdot\|_{+,*}^2\leq
\|\cdot\|_*^2 \end{equation} due to (\ref{wweq2}). Invoking
Proposition \ref{char}, the function $v(\xi)={\kappa\over
2}\|\xi\|_{*,+}^2:B_*\to\bR$ satisfies
$$
\forall (\xi,\eta\in\inter B_*,\;x\in\partial v(\xi),\;y\in\partial
v(\eta)):\; \langle \xi-\eta,x-y\rangle \geq
{\kappa\over\kappa_+}\|\xi-\eta\|_{+,*}^2,
$$
and thus satisfies (\ref{equivalent}.$a$) due to (\ref{notethat}).
At the same time, $$\max\limits_{B_*}v-\min\limits_{B_*}v=
{\kappa\over 2}\max\limits_{\|\xi\|_*\leq1}\|\xi\|_{+,*}^2\leq
{\kappa\over 2},$$ where the concluding inequality is due to
(\ref{notethat}). (i) is proved.
\paragraph{(ii):} Let $v(\cdot)$ satisfy (\ref{equivalent}) and
(\ref{issuchthat}); clearly, the function ${1\over
2}[v(\xi)+v(-\xi)]-v(0)$ also satisfy these relations; thus, we
can assume w.l.o.g. that $v(\xi)=v(-\xi)$ and $v(0)=0$. Let $V$ be
the Legendre transform of $v(\cdot)\bigg|_{B_*}$, that is,
$$
V(x)=\max\limits_{\|\xi\|_*\leq 1}\left[\langle \xi,x\rangle -
v(x)\right].
$$
By the standard properties of the Legendre transform,
(\ref{equivalent}) implies that $V$ is a continuously
differentiable convex  function on $E$ such that
$$
V'(x)=\argmin_{\xi\in
B_*}\left[\langle\xi,x\rangle-v(\xi)\right]\in B_*\ \hbox{and}\
\|V'(x)-V'(y)\|_*\leq \|x-y\|\,\,\forall x,y.
$$
In addition, we clearly have $V(x)=V(-x)$ and $\|x\|-{\kappa\over
2}\leq V(x)\leq \|x\|$ for all $x$ by (\ref{issuchthat}).
Convolving $V$ with a smooth symmetric w.r.t. the origin
nonnegative kernel with unit integral and small support and
subtracting a constant to make function vanish at the origin, we
see that for every $\epsilon>0$ there exists a C$^\infty$ convex
function $W=W_\epsilon$ on $E$ such that for all $x\in E$ one has
\begin{equation}\label{suchthat2}
\begin{array}{ll}
(a)&W_\epsilon(x)=W_\epsilon(-x), \,\, W_\epsilon(0)=0;\\
(b)&\|x\|-{\kappa\over 2}-\epsilon \leq W_\epsilon(x)\leq \|x\|+\epsilon\\
(c)&\|W'(x)\|_*\leq 1\\
(d)&0\leq \langle W''(x)dx,dx\rangle \leq \|dx\|^2\,\,\forall
dx\in E.\\
\end{array}
\end{equation}
Assuming $\epsilon\leq \kappa/10$, let us set
$B=\{x:W(x)\leq\kappa\}$. Then $B$ is a closed convex set
symmetric w.r.t. the origin and such that
\begin{equation}\label{suchthat3}
\{x:\|x\|\leq {9\over 10}\kappa\}\subset B\subset \{x:\|x\|\leq
{5\over 2}\kappa\}
\end{equation}
due to (\ref{suchthat2}.$b$). $B$ is the unit ball of certain norm
$r(x)$ on $E$; by (\ref{suchthat3}) we have
\begin{equation}\label{eq16}
{2\over 5} \|x\|\leq \kappa r(x)\leq {10\over 9}\|x\|.
\end{equation}
Setting $L(x)=p^2(x)$, observe that the function $L$ is given by
the equation
$$
V(x/\sqrt{L(x)})=\kappa.
$$
It follows immediately from the Implicit Function Theorem  that $L$
is C$^\infty$ outside of the origin, and since this function is
the square of a norm, it is therefore C$^1$ on the entire space.
Let us compute the second order differential of $L$ at a point
$x\neq0$. Differentiating twice the equation specifying $L$, we
get
\bse
DL(x)[dx]&=&2L{\langle W',dx\rangle\over\langle
W',x\rangle},\\
D^2L(x)[dx,dx]&=&2L\left[{\langle W',dx\rangle\over\langle
W',x\rangle}\right]^2+{2L^{1/2}\over\langle W',x\rangle}\langle
W''\left[dx-{\langle W',dx\rangle\over\langle
W',x\rangle}x\right],\left[dx-{\langle W',dx\rangle\over\langle
W',x\rangle}x\right]\rangle,\\
\hbox{where}&&L=L(x),W'=W'(L^{-1/2}x), W''=W''(L^{-1/2}x).\\
\ese
We claim that \begin{equation}\label{weclaim} x\neq 0\Rightarrow
0\leq D^2L(x)[dx,dx]\leq {27\over\kappa}\|dx\|^2.
\end{equation}
Indeed, $D^2L(x)[dx,dx]$ is homogeneous of degree 0 in $x$, so
that it suffices to verify the required relation when $L(x)=1$,
i.e., when $W(x)=\kappa$. In this case, the required bound is
readily given by the expression for $D^2L$ combined with
(\ref{suchthat2}.$c,d$) and the following observations: (1) for
$x$ in question, we have $\langle W',x\rangle \geq
W(x)-W(0)=\kappa$, and (2) $\|x\|\leq {5\over 2}\kappa$ by
(\ref{suchthat3}).
\par
Setting $\|x\|_+={5\over 2}\kappa r(x)$ and invoking
(\ref{suchthat3}), we have
\begin{equation}\label{wehave12}
\|\cdot\|^2\leq \|\cdot\|_+^2\leq O(1)\|\cdot\|^2,
\end{equation}
while from (\ref{weclaim}) it follows that the function
$f(x)=\|x\|_+^2$ satisfies
$$
\|f'(x)-f'(y)\|_*\leq O(1)\kappa\|x-y\|,
$$
which combines with (\ref{wehave12}) to imply that
$$
\|f'(x)-f'(y)\|_{+,*}\leq O(1)\kappa\|x-y\|_*.
$$
Thus, $(E,\|\cdot\|)$ is $O(1)\kappa$-smooth, whence, by
(\ref{wehave12}), $(E,\|\cdot\|)$ is $O(1)\kappa$-regular. \myqed
\subsubsection{Proof of Proposition \ref{prop28}}
The fact that a subspace of a $\kappa$-smooth/regular space
equipped with the induced norm is $\kappa$-smooth/regular is
evident. As about the factor-space $F=E/L$, note that the space
dual to $(F,\|\cdot\|_F)$ is nothing but the subspace
$L^\perp=\{\xi:\langle\xi,x\rangle=0\,\forall x \in L\}$ in $E^*$
equipped by the norm induced by $\|\cdot\|_*$. Now assume that
$(E,\|\cdot\|_E)$ is $\kappa$-smooth.  By Proposition \ref{char},
it follows that $\|\cdot\|_*$ possesses property (iv) and
therefore its restriction on $L^\perp$ possesses the same
property. Applying Proposition \ref{char} again, we conclude that
$(F,\|\cdot\|_F)$ is $\kappa$-smooth. We see that passing to a
factor-space preserves $\kappa$-smoothness, and since this
transformation preserves also relations like (\ref{wweq2}), it
preserves $\kappa$-regularity as well. \myqed
\subsubsection{Proof of Proposition \ref{calculus}}
\paragraph{(i):} To prove (i), let $p_i(x^i)=\|x^i\|_i^2$.
\subparagraph{A.} Let $\rho\in[2,\infty)$ be such that $\rho\leq
p$, and let $r=\rho/2$. Our local goal is to prove
\begin{lemma}\label{local} The norm $\|\cdot\|$ on $E=E_1\times ...\times E_m$ defined as
$$\|(x^1,...,x^m)\|=\|(\|x^1\|_1,...,\|x^m\|_m)\|_\rho$$ is
$\kappa_+$-smooth, with
\begin{equation} \label{eqnewkappa}
\kappa_+=\kappa+\rho-2
\end{equation}
\end{lemma}
{\bf Proof.} We have
$$
p(x^1,...,x^m)\equiv
\|(\|x^1\|_1,...,\|x^m\|_m)\|_\rho^2=\|(p_1(x^1),...,p_m(x^m))\|_r.
$$
From this observation it immediately follows that $p(\cdot)$ is
continuously differentiable. Indeed, $\rho\geq2$, whence $r\geq1$,
so that the function $\|y\|_r$ is continuously differentiable
everywhere on $\bR^m_+$ except for the origin; the functions
$p_i(x^i)$  are continuously differentiable by assumption.
Consequently, $p(x)$ is
 continuously differentiable everywhere on
$E=E_1\times...\times E_m$, except, perhaps, the origin; the fact
that $p'$ is continuous at the origin is evident.
\par
Invoking Proposition \ref{char}, in order to prove Lemma
\ref{local} it suffices to verify that
\begin{equation}\label{Lip}
\|p'(x)-p'(y)\|_*\leq 2\kappa_+\|x-y\|
\end{equation}
for all $x,y$. Since $p'$ is continuous, it suffices to prove this
relation for a dense in $E\times E$ set of pairs $x,y$, for
example, those for which all blocks $x^i\in E_i$ in $x$ are
nonzero. With such $x$, the segment $[x,y]$ contains finitely many
points $u$ such that at least one of the blocks $u^i$ is zero;
these points split $[x,y]$ into finitely many consecutive
segments, and it suffices to prove that $$\|p'(x')-p'(y')\|_*\leq
2\kappa_+ \|x'-y'\|$$ when $x',y'$ are endpoints of such a segment.
Since $p'$ is continuous, to prove the latter statement is the
same as to prove similar statement for the case when $x',y'$ are
interior points of the segment. The bottom line is as follows: in
order to prove (\ref{Lip}) for all pairs $x,y$, it suffices to
prove the same statement for those pairs $x,y$ for which every
segment $[x^i,y^i]$ does not pass through the origin of the
corresponding $E_i$. \par Let $x,y$ be such that $[x^i,y^i]$ does
not pass through the origin of $E_i$, $i=1,...,m$. Same as in the
item ``(i)$\Rightarrow$(iii)'' of the proof of Proposition
\ref{char}, for every $i$ there exists a sequence of C$^\infty$
convex functions $\{p_i^t(\cdot)>0\}_{t=1}^\infty$ on $E_i$
converging to $p_i(\cdot)$ along with first order derivatives
uniformly on compact sets and such that
\begin{equation}\label{smooth}
|D^2p_i^t(u^i)[h^i,h^i]| \leq 2\kappa\|h^i\|_i^2\,\,\forall
(u^i,h^i\in E_i).
\end{equation}
Functions $p^t(u)=\|(p_1^t(u^1),...,p_m^t(u^m))\|_r$ clearly are
convex, C$^\infty$ (recall that $p_i^t(\cdot)>0$) and converge to
$p(\cdot)$, along with their first order derivatives, uniformly on
compact sets. It follows that
\begin{equation}\label{limit}
\langle p'(y)-p'(x),h\rangle
=\lim\limits_{t\to\infty}\int\limits_0^1D^2p^t(x+t(y-x))[y-x,h]dt.
\end{equation}
 Setting
$F(y_1,...,y_m)=y_1^r+...+y_m^r$, $y\geq0$, we have
$p^t(u)=F^{{1\over r}}(p_1^t(u^1),...,p_m^t(u^m))$. Now let
$u\in[x,y]$, and let $v\in E$. We have
$$
\begin{array}{rcl}
Dp^t(u)[v]&=&r^{-1}F^{{1\over
r}-1}(p_1^t(u^1),...,p_m^t(u^m))\left(\sum\limits_i
r(p_i^t(u^i))^{r-1}Dp_i^t(u^i)[v^i]\right)\\
\Rightarrow D^2p^t(u)[v,v]&=&{1\over r}\underbrace{\left({1\over
r}-1\right)}_{\leq0}F^{{1\over
r}-2}(p_1^t(u^1),...,p_m^t(u^m))\left(\sum\limits_i
r(p_i^t(u^i))^{r-1}Dp_i^t(u^i)[v^i]\right)^2\\
\multicolumn{3}{l}{+F^{{1\over
r}-1}(p_1^t(u^1),...,p_m^t(u^m))\sum\limits_i\left[(r-1)(p_i^t(u^i))^{r-2}
(Dp_i^t(u^i)[v_i])^2+
(p_i^t(u^i))^{r-1}D^2p_i^t(u^i)[v^i,v^i]\right]}\\
\multicolumn{3}{l}{\leq F^{{1\over
r}-1}(p_1^t(u^1),...,p_m^t(u^m))\sum\limits_i\left[(r-1)(p_i^t(u^i))^{r-2}
(Dp_i^t(u^i)[v_i])^2+2\kappa(p_i^t(u^i))^{r-1} p_i(v^i)\right]}\\
\end{array}
$$
whence
\begin{equation}\label{bottom1}
\begin{array}{l}
0\leq D^2p^t(u)[v,v]\\
\leq F^{{1\over
r}-1}(p_1^t(u^1),...,p_m^t(u^m))\sum\limits_i\left[(r-1)(p_i^t(u^i))^{r-2}
(Dp_i^t(u^i)[v_i])^2+2\kappa(p_i^t(u^i))^{r-1} p_i(v^i)\right].\\
\end{array}
\end{equation}
 Taking into account that $p_i(\cdot)$ are bounded away from
zero on $[x,y]$ and that $p_i^t(\cdot)$ converge, along with first
order derivatives, to $p_i(\cdot)$ uniformly on compact sets as
$t\to\infty$, the right hand side in bound (\ref{bottom1})
converges, as $t\to\infty$, uniformly in $u\in[x,y]$ and $v$,
$\|v\|\leq1$, to
$$
\Psi(u,v)= \left(\sum\limits_i
\|u^i\|_i^\rho\right)^{{2\over\rho}-1}\sum\limits_i\left[(r-1)\|u^i\|_i^{\rho-4}
(Dp_i(u^i)[v_i])^2+2\kappa\|u^i\|_i^{\rho-2}\|v^i\|_i^2\right].
$$
By evident reasons, $|Dp_i(u^i)[v_i]|\leq 2\|u^i\|\|v^i\|$, whence
\begin{equation}\label{bottom2}
\begin{array}{rcl}
\Psi(u,v)&\leq&\left(\sum\limits_i
\|u^i\|_i^\rho\right)^{{2\over\rho}-1}\sum\limits_i\left[4(r-1)\|u^i\|_i^{\rho-2}
\|v_i\|_i^2+2\kappa\|u^i\|_i^{\rho-2}\|v^i\|_i^2\right]\\ &=&
\underbrace{[2\rho+2\kappa-4]}_{2\kappa_+}\left(\sum\limits_i
\|u^i\|_i^\rho\right)^{{2\over\rho}-1}\sum\limits_i\|u^i\|_i^{\rho-2}
\|v^i\|_i^2\\
\end{array}
\end{equation}
When $\rho>2$, we have
$$
\begin{array}{rcl}
\sum\limits_i\|u^i\|_i^{\rho-2} \|v^i\|_i^2&\leq&
\left(\sum\limits_i(\|u^i\|_i^{\rho-2})^{{\rho\over\rho-2}}\right)^{{\rho-2\over
\rho}}\left(\sum\limits_i (\|v^i\|_i^2)^{{\rho\over
2}}\right)^{{2\over\rho}}\\
&=&\left(\sum\limits_i\|u^i\|_i^\rho\right)^{{\rho-2\over
\rho}}\left(\sum\limits_i \|v^i\|_i^\rho\right)^{{2\over\rho}},\\
\end{array}
$$
and (\ref{bottom2}) implies that $ \Psi(u,v)\leq 2\kappa_+\|v\|^2.
$ This inequality clearly is valid for $\rho=2$ as well. Recalling
the origin of $\Psi(\cdot,\cdot)$, we conclude that for every
$\epsilon>0$ there exists $t_\epsilon$ such that
$$
t\geq t_\epsilon,u\in[x,y],\|v\|\leq1\Rightarrow 0\leq
D^2p^t(u)[v,v]\leq 2\kappa_+\|v\|^2+\epsilon.
$$
The resulting inequality via the same reasoning as in the proof of
item ``(i)$\Rightarrow$(iii)'' of Proposition \ref{char} implies
that
$$
t\geq t_\epsilon, u\in[x,y]\Rightarrow |D^2p^t(u)[v,w]|\leq
(2\kappa_++\epsilon)\|v\|\|w\|\,\,\forall v,w.
$$
In view of this bound and (\ref{limit}), we conclude that
$$
\langle p'(y)-p'(x),h\rangle\leq(2\kappa_++\epsilon)\|y-x\|\|h\|
$$
for all $h$, whence $\|p'(y)-p'(x)\|_*\leq
(2\kappa_++\epsilon)\|y-x\|$. Since $\epsilon>0$ is arbitrary, we
arrive at (\ref{Lip}). \subparagraph{B.} When $\rho\leq p$, we
have
$$
\|(\|x^1\|_1,...,\|x^m\|_m)\|_p^2\leq
\|(\|x^1\|_1,...,\|x^m\|_m)\|_\rho^2\leq m^{{2\over\rho}-{2\over
p}}\|(\|x^1\|_1,...,\|x^m\|_m)\|_p^2,
$$
which combines with Lemma \ref{local} to imply that the norm in
(i) is $\kappa$-regular with
$\kappa=[\rho+\kappa-2]m^{{2\over\rho}-{2\over p}}$, for every
$\rho\in[2,p]$, and (i) follows.
\paragraph{(ii):} To prove (ii), consider the norm
$|(x^1,...,x^m)|=m^{1/2}\sqrt{\|x^1\|_1^2+...+\|x^m\|_m^2}$ on
$E\times E\times...\times E$. As it is immediately seen, this norm
is $\kappa$-smooth. If, further,
$\|(x^1,...,x^m)\|_\dag=\sum\limits_i\|x^i\|_i$, then
$$
\|x\|_\dag^2\leq |x|^2\leq m\|x\|_\dag^2\,\,\,\,\,\forall x\in
E\times...\times E,
$$
whence $\|\cdot\|_\dag$ is $m\kappa$-regular. The norm in (ii) is
nothing but the restriction of $\|\cdot\|_\dag$ on the image of
$E$ under the embedding $x\mapsto(x,...,x)$ of $E$ into
$E\times...\times E$, and it remains to use Proposition
\ref{prop28}.   \blacksquare
\subsubsection{Proof of Proposition \ref{calculusI}}
\paragraph{A useful lemma} We start with the following fact:
\begin{lemma}\label{lemfact} Let $(E,\|\cdot\|)$ be a
finite-dimensional $\kappa$-regular space. Then there exists
$\kappa$-smooth norm $\|\cdot\|_+$ on $E$ such that
\begin{equation}\label{???}
\forall (x\in E): \|x\|^2\leq \|x\|_+^2\leq 2\|x\|^2.
\end{equation}
\end{lemma}
{\bf Proof.} By definition, there exists $\kappa_+\in[1,\kappa]$
and a norm $\pi(\cdot)$ on $E$ which is $\kappa_+$-smooth and such
that
$$
\forall(x\in E):\|x\|^2\leq \pi^2(x)\leq
\mu\|x\|^2,\,\,\mu=\kappa/\kappa_+,
$$
or, which is the same, \begin{equation}\label{bbeq0} \forall
\xi\in E^*: \mu\pi_*^2(\xi)\geq \|\xi\|_*^2\geq
\pi_*^2(\xi),
\end{equation}
where $E^*$ is the space dual to $E$ and $\pi_*$, $\|\cdot\|_*$
are the norms on $E^*$ conjugate to $\pi$, $\|\cdot\|$,
respectively.
\par
In the case of $\mu\leq 2$, let us take $\|\cdot\|_+\equiv
\pi(\cdot)$, thus getting a $\kappa_+$-smooth (and thus --
$\kappa$-smooth as well) norm on $E$ satisfying (\ref{???}). Now
let $\mu>2$, so that $\gamma=1/(\mu-1)\in(0,1)$. Let us set
$q_*(\xi)=\sqrt{\gamma \aic{}{\mu} \pi_*^2(\xi)+(1-\gamma)\|\xi\|_*^2}$, so
that $q_*(\cdot)$ is a norm on $E^*$. We have
\begin{equation}\label{bbeq1}
\forall \xi\in E^*:q_*^2(\xi)\geq \|\xi\|_*^2\geq{1\over
\gamma\mu+1-\gamma}q_*^2(\xi)={1\over 2}q_*^2(\xi).
\end{equation}
Further, by Proposition \ref{char} we have
$$
\forall (\xi,\eta\in
E^*,\;x\in\partial\pi_*^2(\xi)):\;\pi_*^2(\xi+\eta)\geq\pi_*^2(\xi)+\langle\eta,x\rangle
+{1\over \kappa_+}\pi_*^2(\eta),
$$
whence, due to $\|\xi+\eta\|_*^2\geq \|\xi\|_*^2
+\langle\eta,y\rangle$ for all $\xi,\eta$ and every {\color{blue} $x\in\partial\pi_*^2(\xi)$ and} $y$ from the
subdifferential $D(\xi)$ of $\|\cdot\|_*^2$ at the point $\xi$,
$$
q_*^2(\xi+\eta)\geq q_*^2(\xi)+\langle \eta,\aic{}{\mu\gamma}x+\aic{}{(1-\gamma)}y\rangle
+{\aic{\gamma}{\mu\gamma}\over \kappa_+}\pi_*^2(\eta)\geq q_*^2(\xi)+\langle \eta,\aic{}{\mu\gamma}x+\aic{}{(1-\gamma)}y\rangle +{\gamma\over \kappa_+}q_*^2(\eta)
$$
(note that ${\color{red}\pi_*(\cdot)\geq q_*(\cdot)/\mu}{\color{blue}\pi_*^2(\cdot)\geq q_*^2(\cdot)/\mu}$ by (\ref{bbeq0})). Since
\[\aic{}{\mu\gamma}\partial \pi_*^2(\xi)+\aic{}{(1-\gamma)}D(\xi)=\partial q_*^2(\xi)\] and
${\gamma\over\kappa_+}={1\over (\mu-1)\kappa_+}\geq{1\over
\kappa}$, we get $$ \forall (\xi,\eta\in E^*,z\in\partial
q_*^2(\xi)): q_*^2(\xi+\eta)\geq q_*^2(\xi)+\langle \eta,z\rangle
+{1\over \kappa}q_*^2(\eta). $$  By the same Proposition
\ref{char}, it follows that the norm $\|\cdot\|_+\equiv q(\cdot)$
on $E$ such that $q_*(\cdot)$ is the conjugate of $q(\cdot)$ is
$\kappa$-smooth. At the same time, (\ref{bbeq1}) implies
(\ref{???}). \blacksquare
\paragraph{Proof of Proposition \ref{calculusI}}
is readily given by Lemma \ref{lemfact} combined with the
corresponding items of Proposition \ref{calculus}. E.g., to prove
(i), note that by Lemma \ref{lemfact} we can find $\kappa$-smooth
norms $q_i(\cdot)$ on $E_i$ such that $q_i^2(x^i)\leq
\|x^i\|_i^2\leq 2q_i^2(x^i)$ for every $i$ and all $x^i\in E_i$.
Applying Proposition \ref{calculus}.(i) to the spaces
$(E_i,q_i(\cdot))$, we get that the norm $q(x^1,...,x^m)=
\left(\sum\limits_{i=1}^m q_i^p(x^i)\right)^{1/p}$  on
$E_1\times...\times E_m$ is $\kappa^+$-regular with $\kappa^+$
given by (\ref{kappa_plus}). Taking into account the evident
relation
$$
q^2(x^1,...,x^m)\leq \|(x^1,...,x^m)\|^2\leq 2q^2(x^1,...,x^m)
$$
and recalling the definition of regularity, we conclude that
$\|\cdot\|$ is $\kappa^{++}$-regular, as required. \blacksquare
\subsection{Proof of Theorem \ref{mainthe}}
\subsubsection{Reduction to the case of a smooth norm}
We intend to reduce the situation to the one where $(E,\|\cdot\|)$
is $\kappa$-smooth rather than $\kappa$-regular. Specifically, we
are about to prove the following fact:
\begin{theorem}
\label{sumthe21} Let $(E,\|\cdot\|)$ be $\kappa$-smooth, let
$E$-valued martingale-difference $\xi^\infty$ satisfy
$(\cC_\alpha[\sigma^\infty])$, and let $S_N=\sum_{i=1}^N \xi_i$,
$\sigma^N=[\sigma_1;...;\sigma_N]$. Then
\par {\rm (i)} When $1\leq\alpha\leq2$, one has for all $N\ge 1$ and $\g\ge 0$:
\begin{equation}\label{fineq3}
 \Prob\left\{\|S_N\|\geq
\left[\sqrt{\exp\{1\}\kappa}+\gamma\right]\sqrt{\sum_{i=1}^N\sigma_i^2}\right\}\leq
2\exp\left\{-{1\over
64}\min\left[\gamma^2;\gamma_*^{2-\alpha}\gamma^\alpha\right]\right\},
\end{equation}
where
\begin{equation}\label{fineq1}
\begin{array}{rcl}
\gamma_*\equiv\gamma_*(\alpha,\sigma^N)&=&\left\{\begin{array}{ll}
\begin{array}{l}32\left[{8\alpha_*\over2^{\alpha_*}}\right]^{{\alpha-1\over2-\alpha}}
\left[{\|\sigma^N\|_2\over
\|\sigma^N\|_{\alpha_*}}\right]^{{\alpha\over
2-\alpha}}\geq16\left[{\|\sigma^N\|_2\over
\|\sigma^N\|_{\alpha_*}}\right]^{{\alpha\over 2-\alpha}}\geq16,\\
\alpha_*={\alpha\over \alpha-1},\\
\end{array}
&1<\alpha<2,\\
\lim_{\alpha \to 1+0}\gamma_*(\alpha,\sigma^N)=16{\|\sigma^N\|_2\over\|\sigma^N\|_\infty},&\alpha=1,\\
\lim_{\alpha\to 2-0}\gamma_*(\alpha,\sigma^N)=+\infty,&\alpha=2.\\
\end{array}\right.\\
\end{array}
\end{equation}
\par {\rm (ii)} When $\alpha=2$, the bound {\rm (\ref{fineq3})}
improves to
\begin{equation}\label{fineq3_impr}
(\forall N\geq1,\gamma\geq0): \Prob\left\{\|S_N\|\geq
\left[\sqrt{\kappa}+\gamma\right]\sqrt{\sum_{i=1}^N\sigma_i^2}\right\}\leq
\exp\{-\gamma^2/3\}.
\end{equation}
\par
{\rm (iii)} When the condition
$\bE_{i-1}\left\{\exp\{\|\xi_i\|^2\sigma_i^{-2}\}\right\}\leq\exp\{1\}$
in $(\cC_2[\sigma^\infty])$ is strengthened to
$\|\xi_i\|\leq\sigma_i$ almost surely, $i=1,2,...,$ the bound {\rm
(\ref{fineq3})} improves to
\begin{equation}\label{fineq3_impr_impr}
(\forall N\geq1,\gamma\geq0): \Prob\left\{\|S_N\|\geq
\left[\sqrt{\kappa}+\gamma\right]\sqrt{\sum_{i=1}^N\sigma_i^2}\right\}\leq
\exp\{-\gamma^2/2\}.
\end{equation}
\end{theorem}
It is immediately seen that Theorem \ref{sumthe21} implies Theorem
\ref{mainthe}. Indeed, if $(E,\|\cdot\|)$ is $\kappa$-regular, by
Lemma \ref{lemfact} there exists a  norm $\|\cdot\|_+$ on $E$ such
that $(E,\|\cdot\|_+)$ is $\kappa$-smooth and (\ref{???}) holds
true. Setting $\widehat{\sigma}_i=\sqrt{2}\sigma_i$, observe that
(\ref{???}) combines with $(\cC_\alpha[\sigma^\infty])$ to imply
that
$\bE_{i-1}\left\{\exp\{\|\xi_i\|_+^2\widehat{\sigma}_i^{-2}\}\right\}\leq
\exp\{1\}$. Applying Theorem \ref{sumthe21}.(i) to the
$\kappa$-smooth space $(E,\|\cdot\|_+)$ and $\widehat{\sigma}_i$
in the role of $\sigma_i$ and taking into account that
$\|S_N\|\leq\|S_N\|_+$, we see that Theorem \ref{mainthe}.(i) is
an immediate corollary of Theorem \ref{sumthe21}.(i), and
similarly for Theorem \ref{mainthe}.(ii-iii). \subsubsection{Proof
of Theorem \ref{sumthe21}: preliminaries} In the sequel, we need
the following (essentially, well-known) fact.
\begin{proposition}\label{boundprob} Let $\psi_i$, $i=1,...,N$, be Borel functions
on $\Omega$ such that $\psi_i$ is $\cF_i$-measurable,  let
$\alpha\in [1,2]$, and let $\mu_i$, $\nu_i>0$ be deterministic
reals. Assume that almost surely one has
\begin{equation}\label{geq1}
\bE_{i-1}\{\psi_i\}\leq\mu_i,
\bE_{i-1}\left\{\exp\{|\psi_i|^\alpha/\nu_i^\alpha\}\right\}\leq
\exp\{1\},\,1\leq i\leq N.
\end{equation}
Then for every $\gamma\geq 0$ one has
\begin{equation}\label{Q.E.D}
\Prob\left\{\sum_{i=1}^N\psi_i>\sum_{i=1}^N\mu_i+\gamma\sqrt{\sum_{i=1}^N\nu_i^2}\right\}
\leq
2\exp\{-{1\over64}\min\left[\gamma^2,\gamma_*^{2-\alpha}\gamma^\alpha\right]\},
\end{equation}
where
\begin{equation}\label{Q.E.D.1}
\gamma_*\equiv\gamma_*(\alpha,\nu^N)=\left\{\begin{array}{ll}
\begin{array}{l}32\left[{8\alpha_*\over
2^{\alpha_*}}\right]^{{\alpha-1\over2-\alpha}}\left[{\|\nu^N\|_2\over\|\nu^N\|_{\alpha_*}}
\right]^{{\alpha\over2-\alpha}} \geq
16\left[{\|\nu^N\|_2\over\|\nu^N\|_{\alpha_*}}
\right]^{{\alpha\over2-\alpha}}\geq16,\\
\left[\alpha_*={\alpha\over\alpha-1},\nu^N=[\nu_1;...;\nu_N]\right],\\
\end{array}&1<\alpha<2,\\
\lim_{\alpha \to
1+0}\gamma_*(\alpha,\nu^N)=16{\|\nu^N\|_2\over\|\nu^N\|_\infty},&\alpha=1,\\
\lim_{\alpha\to 2-0}\gamma_*(\alpha,\nu^N)=+\infty,&\alpha=2.\\
\end{array}\right.
\end{equation}
\end{proposition}
To make the text self-contained, here is the proof.
\paragraph{0$^0$.} Till item {\bf 4$^0$} of the proof, we restrict ourselves with the case when
$1<\alpha<2$. Besides this, by evident homogeneity reasons we may
assume w.l.o.g. that $\nu\equiv\sum_{i=1}^N\nu_i^2=1$.
\paragraph{1$^0$.} We
start with the following
\begin{lemma}\label{ilemma} Let $\alpha\in(1,2)$, $\nu>0$ and  $\psi$ be a real-valued
random variable such that \begin{equation}\label{ieq5}
\bE\{\exp\{|\psi/\nu|^\alpha\}\}\leq\exp\{1\}.
\end{equation}
Then
\begin{equation}\label{ieq6}
t\geq 0\Rightarrow \ln\left(\bE\{\exp\{t\psi\}\}\right)\leq
t\bE\{\psi\}+8
(t\nu)^2+2^{\alpha_*}\alpha_*^{-1}|t\nu|^{\alpha_*},\quad
\alpha_*={\alpha\over\alpha-1}.
\end{equation}
\end{lemma}
{\bf Proof.} {\bf 1)} Let $t\geq0$ be fixed. W.l.o.g. we can
assume that $\nu=1$. By Young inequality, we have
$$ t\psi=(2 t)(\psi/2)\leq {|\psi/2|^\alpha\over\alpha}+{(2
t)^{\alpha_*}\over\alpha_*};$$ since $\alpha^{-1}(1/2)^\alpha<1
$
and $\nu=1$, we have
$\bE\{\exp\{\alpha^{-1}|\psi/2|^\alpha\}\}\leq\exp\{\alpha^{-1}(1/2)^\alpha\}$,
whence
\bse
\bE\{\exp\{t\psi\}\}\leq\bE\{\exp\{\alpha^{-1}|\psi/2|^\alpha
+\alpha_*^{-1}(2 t)^{\alpha_*}\}\}\leq
\exp\{\alpha^{-1}(1/2)^\alpha+\alpha_*^{-1}(2t)^{\alpha_*}\}.
\ese
{\bf 2)}  Let $f(t)=\bE\{\exp\{t\psi\}\}$. Since $\alpha>1$, $f$
is a C$^\infty$ function on the axis such that $f(0)=1$,
$f'(0)=\bE\{\psi\}$ and
\bse
f''(t)=\bE\left\{\exp\{t\psi\}\psi^2\right\}
\ese
It is easily seen that
\bse
0\leq t\leq 1/4\Rightarrow \exp\{t
|s|\}s^2\leq\exp\{|s|^\alpha\}\,\forall s,
\ese
whence under the premise of Lemma \ref{ilemma} one has
$$
0\leq t\leq 1/4\Rightarrow f''(t)\leq \exp\{1\}
$$
(recall that $\nu=1$). It follows that $$0\leq t\leq
1/4\Rightarrow f(t)\leq 1+t\bE\{\psi\}+{\exp\{1\}\over
2}t^2\leq\exp\{t\bE\{\psi\}+{\exp\{1\}\over 2}t^2\}.
$$
Thus, one has
\begin{equation}\label{lrcl}
\begin{array}{lrcl}
(a)&0\leq t\leq 1/4\Rightarrow \ln f(t)&\leq&
t\bE\{\psi\}+{\exp\{1\}\over 2}t^2,\\
(b)&t\geq 0\Rightarrow \ln f(t)&\leq&
\alpha^{-1}(1/2)^\alpha+\alpha_*^{-1}(2t)^{\alpha_*}.\\
\end{array}
\end{equation}
Since $8 t^2\geq {\exp\{1\}\over2}t^2$ and $8 t^2\geq
\alpha^{-1}(1/2)^\alpha$ when $t\geq1/4$, (\ref{lrcl}) implies
 (\ref{ieq6}).
\paragraph{2$^0$.} Since $\alpha>1$, we have for all $t\geq0$
$$
\begin{array}{l}
\bE\left\{\exp\{t\sum_{i=1}^n\psi_j\}\right\}=
\bE\left\{\exp\{t\sum_{i=1}^{n-1}\psi_j\}\bE_{n-1}\{\exp\{t\psi_n\}\}\right\}\\
\leq
\bE\left\{\exp\{\exp\{t\sum_{i=1}^{n-1}\psi_j\}\right\}\exp\{\mu_nt+8
(t\nu_n)^2+\alpha_*^{-1}2^{\alpha_*}(t\nu_n)^{\alpha_*}\},\\
\end{array}
$$
whence
\bse
\begin{array}{l}
\ln\left(\bE\{t\sum_{i=1}^N\psi_i\}\right\}\leq
A_Nt+B_Nt^2+C_Nt^{{\alpha\over\alpha-1}},\\
A_N=\sum_{i=1}^N\mu_i,\, B_N=8\sum_{i=1}^N\nu_i^2,\,
C_N=\alpha_*^{-1}2^{\alpha_*}\sum_{i=1}^N\nu_i^{\alpha_*}.\\
\end{array}
\ese

\paragraph{3$^0$.} Recall that we are in the situation $\sum_{i=1}^N\nu_i^2=1$.
We have for all $t>0$:
\bse
\Prob\left\{\Psi_N>A_N+\gamma\nu\right\}&=&\Prob\left\{\exp\{t\Psi_N\}
>\exp\{tA_N+t\gamma\}\right\}\\
&\leq&
\bE\left\{\exp\{t\Psi_N\}\right\}\exp\{-tA_N-t\gamma\}\leq
\exp\{B_Nt^2+C_Nt^{{\alpha\over\alpha-1}}-t\gamma\},
\ese
whence
$$
\Prob\{\Psi_N>A_N+\gamma\}\leq\inf_{t>0}\exp\{B_Nt^2+C_Nt^{{\alpha\over\alpha-1}}-t\gamma\}.
$$
whence also
\bse
\ln\left(\Prob\{\Psi_N>A_N+\gamma\}\right)\leq
\ln(2)+\inf_{t>0}\big[
\underbrace{\max[2B_Nt^2,2C_Nt^{{\alpha\over\alpha-1}}]}_{\phi(t)}-\gamma
t\big]\equiv \ln(2)-\phi_*(\gamma),
\ese
where $\phi_*$ is the Legendre transform of $\phi$,
$\Dom\phi=[0,\infty)$. Let $t_*=t_*(\alpha)$ be the unique
positive root of the equation $B_Nt^2=C_Nt^{\alpha_*}$, that is,
$$ t_*=(B_N/C_N)^{{\alpha-1\over 2-\alpha}}. $$ The function
$\phi(t)$ is strongly convex on $[0,\infty)$, equals $2B_Nt^2$ to
the left of $t_*$ and equals $2C_Nt^{\alpha_*}$ to the right of
$t_*$. Let $\gamma_-=\gamma_-(\alpha)$ be the left, and
$\gamma_+=\gamma_+(\alpha)$ be the right derivative of $\phi$ at
$t_*$, so that $$
4B_Nt_*=\gamma_-\leq\gamma_+=2C_N\alpha_*t_*^{{1\over\alpha-1}}.
$$
The function $\phi_*(\gamma)$ is as follows: since $\phi$ is
strongly convex on $[0,\infty)$, $\phi'(0)=0$ and
$\phi(t)/t\to\infty$ as $t\to\infty$, $\phi_*$ is continuously
differentiable and convex on $[0,\infty)$; when $0\leq\gamma
\leq\gamma_-$, $\phi_*$ coincides with the Legendre transform
$\phi_{*,-}(\gamma)={1\over 8B_N}\gamma^2$ of the function
$2B_Nt^2$ on the axis; when $\gamma\geq \gamma_+$, $\phi_*$
coincides with the Legendre transform
$\phi_{*,+}(\gamma)={(2C_N)^{1-\alpha}\over\alpha}\gamma^\alpha$
of the function $2C_N|t|^{\alpha_*}$ on the axis. In the segment
$[\gamma_-,\gamma_+]$ $\phi_*$ is linear with the slope
$\phi_{*,-}^\prime(\gamma_-)=\phi_{*,+}^\prime(\gamma_+)=t_*$. Now
let $\theta=\phi_{*,-}(\gamma_-)/\phi_{*,+}(\gamma_-)$, and let
$\omega(\gamma)=\theta\phi_{*,+}(\gamma)$. Observe that
$\omega(\gamma)\leq\phi_*(\gamma)$ when $\gamma\geq\gamma_-$.
\begin{quote}
Indeed, at the point $\gamma_+$ the functions $\phi_{*,+}$ and
$\phi_*$ have equal values and equal derivatives, and since
$\phi_*$ is linear in $\Delta=[\gamma_-,\gamma_+]$, we conclude
from convexity of $\phi_{*,+}(\cdot)$  that
$\phi_{*,+}(\gamma)\geq \phi_*(\gamma)$ on $\Delta$, while $0\leq
\phi_{*,+}^\prime(\gamma)\leq\phi_*^\prime(\gamma)\equiv\phi_{*,+}^\prime(\gamma_+)$
on $\Delta$. Therefore $\theta\leq1$, and since $\phi_*^\prime$ is
nondecreasing, we have $\omega'(\gamma)\leq \phi_*^\prime(\gamma)$
on $\Delta$. Since $\omega(\gamma_-)=\phi_*(\gamma_-)$, we
conclude that $\omega\leq \phi_*$ everywhere on $\Delta$. Since
$\theta<1$ and $\phi_{*,+}$ is positive, when $\gamma\geq\gamma_+$
we have $\omega(\gamma)\leq\phi_{*,+}(\gamma)=\phi_*(\gamma)$.
\end{quote}
The bottom line is that
\bse
\phi_*(\gamma)\geq\left\{\begin{array}{ll} {1\over
8B_N}\gamma^2&,0\leq\gamma\leq\gamma_-\\
D_N\gamma^\alpha&,\gamma\geq\gamma_-\\
\end{array}\right., \, D_N={\phi_{*,-}(\gamma_-)\over
\gamma_-^\alpha}
\ese
Recalling the definition of $A_N$, $B_N$. $C_N$, we arrive at
(\ref{Q.E.D}) -- (\ref{Q.E.D.1}).
\paragraph{4$^0$.} We have proved the assertion of Proposition in
the case of $1<\alpha<2$. This combines with the standard
approximation arguments to yield the assertion in the cases of
$\alpha=1$ and $\alpha=2$.
\subsubsection{Completing the proof of Theorem \ref{sumthe21}}
\paragraph{1$^0$: Preparations.} Given $\kappa$-smooth space $(E,\|\cdot\|)$, let us set
$$
V(\xi)=\left\{\begin{array}{ll}{1\over 2}\|\xi\|^2&,\|\xi\|\leq
1\\
\|\xi\|-{1\over 2}&,\|\xi\|\geq1\\ \end{array}\right.,\quad
V_\beta(\xi)=\beta V(\xi/\beta)\quad[\beta>0],\quad v(x)={1\over
2}\|x\|_*^2.
$$
Observe that \begin{enumerate} \item  $V_{\beta}(\cdot)$ is the
Legendre transform of the restriction of $\beta v(\cdot)$ on the
$\|\cdot\|_*$-unit ball, whence $\|V_\beta^\prime(\xi)\|_*\leq1$
for all $\beta>0$ and all $\xi$, and
\begin{equation}\label{jeq3}
\|x\|_*\leq 1\Rightarrow \langle x,\xi\rangle\leq \beta
v(x)+V_\beta(\xi)\leq {\beta\over 2}+V_\beta(\xi)\,\forall \xi.
\end{equation}
 \item $V(\cdot)$ is continuously differentiable with
$\|V'(\xi)-V'(\eta)\|_*\leq \kappa\|\xi-\eta\|$ and is Lipschitz
continuous, with constant 1, w.r.t. $\|\cdot\|$;
\begin{quote}
{\small The second claim is evident. To prove the first, note that
the function $v(\cdot)$ on the entire $\bR^n$ is strongly convex
w.r.t. $\|\cdot\|_*$ with parameter $\kappa^{-1}$, whence, of
course, so is the function $\hat{v}$  which is equal to $v$ in the
unit ball and is $+\infty$ outside of this ball. Given $\xi,\eta$
and setting $x=V'(\xi)$, $y=V'(\eta)$, we have $\xi\in\partial
\hat{v}(\xi)$, $\eta\in\partial \hat{v}(y)$, whence
$$\|\xi-\eta\|\|x-y\|_*\geq \langle x-y,\xi-\eta\rangle\geq
\kappa^{-1}\|x-y\|_*^2,$$ so that
$$\|V'(\xi)-V'(y)\|_*=\|x-y\|_*\leq \kappa\|\xi-\eta\|.$$}
\end{quote}
\item  One has
\begin{equation}\label{jeq5}
\begin{array}{ll}
(a)&|V_\beta(\xi+\eta)-V_\beta(\xi)|\leq \|\eta\|\\
(b)&V_\beta(\xi+\eta)-V_\beta(\xi)\leq \langle
V_\beta^\prime(\xi),\eta\rangle +{\kappa\over2\beta}\|\eta\|^2.\\
\end{array}
\end{equation}
\begin{quote}
{\small It clearly suffices to consider the case of $\beta=1$,
that is, $V_\beta\equiv V$. By the second claim in item 2, $V$ is
Lipschitz continuous with constant $1$ w.r.t. the norm
$\|\cdot\|$, which implies (\ref{jeq5}.$a$). Relation
(\ref{jeq5}.$b$) is readily given by the Lipschitz continuity of
$V'$, see the first claim in item 2.}
\end{quote}
\end{enumerate}
\paragraph{2$^0$: Proof of Theorem \ref{sumthe21}.(i).}
Let us fix $\beta>0$ and set \bse
S_n=\sum_{i=1}^n\xi_i,\,\,
a_n=V_\beta^\prime(S_{n-1}),\,\,\psi_n=V_\beta(S_n)-V_\beta(S_{n-1}),
\ese
so that $a_n$ is $\cF_{n-1}$-measurable, and $\psi_n$ is
$\cF_n$-measurable. By (\ref{jeq5}.$a$) we have $|\psi_n|\leq
\|\xi_n\|$, whence
\begin{equation}\label{whencea}
\bE_{n-1}\left\{\exp\{|\psi_n|^\alpha/\sigma_n^\alpha\}\right\}\leq
\exp\{1\},
\end{equation}
while by (\ref{jeq5}.$b$) we have
$$
\begin{array}{l}
\bE_{n-1}\left\{\psi_n\right\}\leq\bE_{n-1}\left\{\langle
a_n,\xi_n\rangle +{\kappa\over 2\beta}\|\xi_n\|^2\right\}=
\bE_{n-1}\left\{\langle a_n,\xi_n\rangle +{\kappa\over
2\beta}\|\xi_n\|^2\right\}\\
=\bE_{n-1}\left\{{\kappa\over 2\beta}\|\xi_n\|^2\right\}
\;\;\;{\hbox{[since $a_n$ is $\cF_{n-1}$-measurable
and $\bE_{n-1}\left\{\xi_n\right\}=0$]}}\\
\leq{\kappa\over 2\beta}\sigma_n^2\exp\{1\}.\\
\end{array}
$$
The concluding inequality above can be justified as follows:
setting $\zeta_n=\|\xi_n\|/\sigma_n$, we have
$\bE_{n-1}\left\{\exp\{\zeta_n^\alpha\}\right\}\leq\exp\{1\}$. At
the same time, it is immediately seen that
$$s^2\leq(\alpha\exp\{1\}/2)^{-2/\alpha}\exp\{|s|^\alpha\}$$ for all
$s$, and since $(\alpha\exp\{1\}/2)^{-2/\alpha}\leq1$ when
$1\leq\alpha\leq2$, we get $\bE_{n-1}\{\zeta_n^2\}\leq
\bE_{n-1}\left\{\exp\{|\zeta_n|^\alpha\}\right\}$. Thus, we arrive
at
\begin{equation}\label{whenceb}
\bE_{n-1}\left\{\psi_n\right\}\leq\mu_n:=\exp\{1\}\sigma_n^2.
\end{equation}
\par
Invoking (\ref{jeq3}), we get
$$
\|S_N\|\leq{\beta\over 2}+V_\beta(S_N)={\beta\over
2}+\sum_{i=1}^N\psi_i.
$$
Taking into account (\ref{whencea}), (\ref{whenceb}) and applying
Proposition \ref{boundprob}, we arrive at
$$
\begin{array}{l}
\forall \gamma\geq0: \Prob\left\{\|S_N\|\geq\left[{\beta\over
2}+{\kappa\exp\{1\}\sum_{i=1}^N\sigma_i^2\over
2\beta}\right]+\gamma\sqrt{\sum_{i=1}^N\sigma_i^2}\right\}\leq
2\exp\{-{1\over
64}\min[\gamma^2,\gamma_*^{2-\alpha}\gamma^\alpha]\},
\end{array}
$$
with $\gamma_*=\gamma_*(\alpha,\sigma^N)$  given by
(\ref{Q.E.D.1}). Optimizing this bound in $\beta>0$, we arrive at
 (\ref{fineq3}). Theorem \ref{sumthe21}.(i) is proved.
\paragraph{3$^0$: Proof of Theorem \ref{sumthe21}.(ii-iii).} These results are given by exactly the same reasoning as above, with
the  role of Proposition \ref{boundprob} played by the following
statement:
\begin{proposition}\label{boundprobii} Let $\psi_i$, $i=1,...,N$, be Borel functions
on $\Omega$ such that $\psi_i$ is $\cF_i$-measurable,  and let
$\mu_i\geq0$, $\nu_i>0$ be deterministic reals. Assume that almost
surely one has
\bse
\forall i: \bE_{i-1}\{\psi_i\}\leq\mu_i,
\ese
and either \begin{equation}\label{onehasiia} \forall i:
\bE_{i-1}\left\{\exp\{\psi_i^2/\nu_i^2\}\right\}\leq \exp\{1\},
\end{equation}
or
\begin{equation}\label{onehasiib} \forall i:
|\psi_i|\leq\nu_i.
\end{equation}
 Then for every
$\gamma\geq 0$ one has
\begin{equation}\label{Q.E.Dii}
\Prob\left\{\sum_{i=1}^N\psi_i>\sum_{i=1}^N\mu_i+\gamma\sqrt{\sum_{i=1}^N\nu_i^2}\right\}
\leq \left\{\begin{array}{ll}\exp\{-\gamma^2/3\},&\hbox{\rm case
of (\ref{onehasiia})}\\
\exp\{-\gamma^2/2\},&\hbox{\rm case of (\ref{onehasiib})}\\
\end{array}\right..
\end{equation}
\end{proposition}
{\bf Proof.} Let (\ref{onehasiia}) be the case. It is immediately
seen that $\exp\{s\}\leq s+\exp\{9s^2/16\}$ for all $s$. We
conclude that if $0\leq t\leq {4\over 3\nu_i}$, then
\be\label{weconclude} \bE_{i-1}\left\{ \exp\{t\psi_i\}\right\}&\leq&
t\mu_i+\bE_{i-1}\left\{ \exp\{9t^2\psi^2_i/16\}\right\}\nn
& \leq&
t\mu_i+\exp\{9t^2\nu_i^2/16\}\leq \exp\{t\mu_i +9t^2\nu^2_i/16\}.
\ee{weconclude}
 Besides this, we have $tx\le {3t^2\nu_i^2\over
8}+{2x^2\over 3\nu_i^2}$, so that $$ \bE_{i-1}\left\{
\exp\{t\psi_i\}\right\}\le \exp\left\{{3t^2\nu_i^2\over 8}+{2\over
3}\right\},$$ and the latter quantity is $\leq
\exp({3t^2\nu_i^2\over4})$ when $t\geq {4\over 3\nu_i}$. Invoking
\rf{weconclude}, we arrive at
\begin{equation}\label{where}
t\geq0\Rightarrow\bE_{n-1}\left\{\exp\{t\phi_n\}\right\}\leq\exp\{t\mu_i+3t^2\nu_n^2/4\}.
\end{equation}
It follows that
$$
\begin{array}{c}
\bE \exp\left\{t\sum_{i=1}^n \psi_i\right\}= \bE \left\{ \bE_{n-1}
\left\{\exp\left\{t\sum_{i=1}^n \psi_i\right\}\right\}\right\}
\leq \bE\left\{\exp\left\{t\sum_{i=1}^{n-1} \psi_i\right\}\right\}
\exp(t\mu_n+3t^2\nu_n^2/4),\end{array} $$
 whence $$\begin{array}{c} t\geq0\Rightarrow
\bE\left\{\exp\{t\sum_{i=1}^N\psi_i\}\right\} \leq
\exp\left\{t\sum_{i=1}^N\mu_i+{3t^2\over
4}\sum_{i=1}^N\nu_i^2\right\}.
\end{array}
$$
Therefore for $\gamma\geq0$ we get
$$
\begin{array}{l}
\Prob\left\{\sum_{i=1}^N\psi_i>\sum_{i=1}^N\mu_i+\gamma\sqrt{\sum_{i=1}^N\nu_i^2}\right\}\\
\leq\min_{t>0}\left[\bE\left\{\exp\{t\sum_{i=1}^N\psi_i\}\right\}\exp\{-t\sum_{i=1}^N\mu_i-t\gamma\sqrt{\sum_{i=1}^N\nu_i^2}
\}\right]\\
\leq\min_{t>0} \exp\{t\sum_{i=1}^N\mu_i+{3t^2\over
4}\sum_{i=1}^N\nu_i^2-t\sum_{i=1}^N\mu_i-t\gamma\sqrt{\sum_{i=1}^N\nu_i^2}\}=\exp\{-\gamma^2/3\}\\
\end{array},
$$
as required in the first bound in (\ref{Q.E.Dii}). In the case of
(\ref{onehasiib}), by Azuma-Hoeffding's inequality \cite{azuma}, we have
$$
\forall t\geq0: \bE_{n-1}\left\{\exp\{t\phi_i\}\right\}\leq
\exp\{t\mu_i+\sigma_i^2/2\};
$$
with this relation in the role of (\ref{where}), the above
reasoning results in the second bound in (\ref{Q.E.Dii}). \myqed

\end{document}